\documentclass[a4paper]{article}
\usepackage{amsthm,amssymb,amsmath,enumerate,graphicx,epsf,bbold}

\newcommand{\COLORON}{1}
\newcommand{\NOTESON}{0}
\newcommand{\Debug}{0}
\usepackage[usenames,dvips]{color} 
\usepackage{amsthm,amssymb,amsmath,enumerate,graphicx,epsf,mathbbol,stmaryrd}
\usepackage[dvips, bookmarks, colorlinks=false, breaklinks=true]{hyperref}

\newcommand{\comment}[1]{}
\newcommand{\COMMENT}[1]{}

\definecolor{darkgray}{rgb}{0.3,0.3,0.3}
\newcommand{\defi}[1]{{\color{darkgray}\emph{#1}}}

\newcommand{\acknowledgement}{\section*{Acknowledgement}}

\comment{
	\begin{lemma}\label{}	
\end{lemma}
\begin{proof}

\end{proof}

\begin{theorem}\label{}
\end{theorem} 
\begin{proof} 	

\end{proof}

}

\newtheorem{proposition}{Proposition}[section]
\newtheorem{definition}[proposition]{Definition}
\newtheorem{theorem}[proposition]{Theorem}
\newtheorem{corollary}[proposition]{Corollary}
\newtheorem{lemma}[proposition]{Lemma}
\newtheorem{observation}[proposition]{Observation}
\newtheorem{conjecture}[proposition]{Conjecture}

\newtheorem{problem}[proposition]{Problem}

\newtheorem{examp}[proposition]{Example}

\newtheorem*{noclaim}{Claim}

\newcommand{\example}[2]{\begin{examp} \label{#1} {{#2}}\end{examp}}

\newcommand{\FIG}{0}

\ifnum \NOTESON = 1 \newcommand{\note}[1]{ 

	\ 

	{\color{blue} \hspace*{-60pt} NOTE: \color{Turquoise}{\small  \tt \begin{minipage}[c]{1.1\textwidth}  #1 \end{minipage} \ignorespacesafterend }} 
	
	\ 
	
	}
\else \newcommand{\note}[1]{} \fi

\newcommand{\afsubm}[1]{ \ifnum \Debug = 1 {\mymargin{#1}}
\fi} 

\ifnum \Debug = 1 
\else  \fi

\ifnum \FIG = 1 \newcommand{\fig}[1]{Figure ``{#1}''}
\else \newcommand{\fig}[1]{Figure~\ref{#1}} \fi

\ifnum \Debug = 1 \usepackage[notref,notcite]{showkeys}
\fi

\ifnum \COLORON = 0 \renewcommand{\color}[1]{}
\fi

\newcommand{\showFig}[2]{
   \begin{figure}[htbp]
   \centering
   \noindent
   \epsfbox{#1.eps}
   \caption{\small #2}
   \label{#1}
   \end{figure}
}

\newcommand{\N}{\ensuremath{\mathbb N}}
\newcommand{\R}{\ensuremath{\mathbb R}}
\newcommand{\Z}{\ensuremath{\mathbb Z}}
\newcommand{\Q}{\ensuremath{\mathbb Q}}

\newcommand{\cb}{\ensuremath{\mathcal B}}

\newcommand{\cf}{\ensuremath{\mathcal F}}

\newcommand{\ci}{\ensuremath{\mathcal I}}

\newcommand{\ck}{\ensuremath{\mathcal K}}

\newcommand{\cm}{\ensuremath{\mathcal M}}
\newcommand{\cn}{\ensuremath{\mathcal N}}

\newcommand{\alp}{\ensuremath{\alpha}}
\newcommand{\bet}{\ensuremath{\beta}}
\newcommand{\gam}{\ensuremath{\gamma}}
\newcommand{\del}{\ensuremath{\delta}}
\newcommand{\eps}{\ensuremath{\epsilon}}

\newcommand{\lam}{\ensuremath{\lambda}}
\newcommand{\sig}{\ensuremath{\sigma}}

\newcommand{\CCC}{\ensuremath{\mathcal C}}

\newcommand{\ccg}{\ensuremath{\mathcal C(G)}}
\newcommand{\fcg}{\ensuremath{|G|}}

\newcommand{\zero}{\mathbb 0}
\newcommand{\sm}{\backslash}
\newcommand{\restr}{\upharpoonright}
\newcommand{\noproof}{\qed}

\newcommand{\nin}{\ensuremath{{n\in\N}}}
\newcommand{\iin}{\ensuremath{{i\in\N}}}
\newcommand{\unin}{\ensuremath{[0,1]}}

\newcommand{\sgl}[1]{\ensuremath{\{#1\}}}

\newcommand{\seq}[1]{\ensuremath{(#1_i)_{i\in\N}}} 
\newcommand{\sseq}[2]{\ensuremath{(#1_i)_{i\in #2}}} 
\newcommand{\susq}[2]{\ensuremath{(#1_{#2_i})_{\iin}}} 
\newcommand{\oseq}[2]{\ensuremath{(#1_\alp)_{\alp< #2}}} 
\newcommand{\fml}[1]{\ensuremath{\{#1_i\}_{i\in I}}} 
\newcommand{\ffml}[2]{\ensuremath{\{#1_i\}_{i\in #2}}} 
\newcommand{\ofml}[2]{\ensuremath{\{#1_\alp\}_{\alp< #2}}} 

\newcommand{\g}{\ensuremath{G\ }}
\newcommand{\G}{\ensuremath{G}}

\newcommand{\s}{s}

\newcommand{\ceil}[1]{\ensuremath{\lceil #1 \rceil}}

\newcommand{\Lr}[1]{Lemma~\ref{#1}}
\newcommand{\Tr}[1]{Theorem~\ref{#1}}
\newcommand{\Sr}[1]{Section~\ref{#1}}

\newcommand{\Prb}[1]{Problem~\ref{#1}}
\newcommand{\Cr}[1]{Corollary~\ref{#1}}

\newcommand{\Or}[1]{Observation~\ref{#1}}
\newcommand{\Er}[1]{Example~\ref{#1}}
\newcommand{\Dr}[1]{De\-fi\-nition~\ref{#1}}

\newcommand{\lf}{locally finite}
\newcommand{\lfg}{locally finite graph}

\newcommand{\nlfg}{non-locally-finite graph}

\renewcommand{\iff}{if and only if}
\newcommand{\fe}{for every}
\newcommand{\Fe}{For every}

\newcommand{\st}{such that}
\newcommand{\sot}{so that}
\newcommand{\ti}{there is}
\newcommand{\tho}{there holds}

\newcommand{\obda}{without loss of generality}

\newcommand{\wrt}{with respect to}

\newcommand{\tcs}{topological cycle space}

\newcommand{\labequ}[2]{ \begin{equation} \label{#1} #2 \end{equation} } 
\newcommand{\labtequ}[2]{ \begin{equation} \label{#1} 	\begin{minipage}[c]{0.9\textwidth}  #2 \end{minipage} \ignorespacesafterend \end{equation} }

\newcommand{\mymargin}[1]{
  \marginpar{%
    \begin{minipage}{\marginparwidth}\small%
      \begin{flushleft}%
        {\color{blue}#1}%
      \end{flushleft}%
   \end{minipage}%
  }%
}%

\newcommand{\mySection}[2]{}

\newcommand{\DK}{Diestel and K\"uhn}

\newcommand{\ThmMacLaneC}{\cite{maclane37}}

\newcommand{\LemHeineCa}[1]{
	\begin{lemma}[Heine-Cantor Theorem] \label{#1}
	Let $M$ be a compact metric space, and let $f: M \to N$ be a continuous function, where $N$ is a metric space. Then $f$ is uniformly continuous.	\end{lemma}
}

\newcommand{\etal}{et.\ al.}

\newcommand{\lhom}{\ensuremath{\hat{H}_d}}

\newcommand{\lhomx}{\ensuremath{\lhom(X)}}
\newcommand{\lhomn}{\ensuremath{\hat{H}_{d,n}}}

\newcommand{\homp}{\ensuremath{H_d}}

\newcommand{\hompx}{\ensuremath{\homp(X)}}

\newcommand{\shom}{\ensuremath{H_1}}

\newcommand{\hleq}{\preceq}

\newcommand{\xdim}[1]{$#1$-dimensional}
\newcommand{\ndim}{\xdim{n}}
\newcommand{\adim}{1-dimensional}
\newcommand{\prim}{primitive}

\newcommand{\ec}[1]{\ensuremath{\llbracket #1 \rrbracket}}

\newcommand{\Cs}{Cauchy sequence}

\newcommand{\clex}{circlex} 
\newcommand{\rep}{\sig-representative}

\newcommand{\mts}{metric disc}

\newcommand{\cosp}{constant speed}
\newcommand{\nuho}{null-homotopic}

\newcommand{\lcup}{\ensuremath{\oplus}}
\newcommand{\blcup}{\ensuremath{\bigoplus}}

\newcommand{\vol}{area}
\newcommand{\evol}{excess area}

\newcommand{\arex}{area extension}

\newcommand{\unsu}{unconditionally summable}

\title{Cycle decompositions: from graphs to continua}

\author{Agelos Georgakopoulos\thanks{Supported by GIF grant I-879-124.6/2005 and FWF grant P-19115-N18.}
\medskip \\
 {Technische Universit\"at Graz}\\
  {Steyrergasse 30, 8010}\\
  {Graz, Austria}\\
}

\begin{document}
\maketitle

\begin{abstract}
We generalise a fundamental graph-theoretical fact, stating that every element of the cycle space of a graph is a sum of edge-disjoint cycles, to arbitrary continua. To achieve this we replace graph cycles by topological circles, and replace the cycle space of a graph by a new homology group for continua which is a quotient of the first singular homology group $H_1$. This homology seems to be particularly apt for studying spaces with infinitely generated $H_1$, e.g.\ infinite graphs or fractals. 
\end{abstract}

\section{Introduction}

\subsection{Overview}
In a recent series of papers, Diestel \etal\ showed that many well-known theorems about cycles in finite graphs remain true for infinite graphs provided one replaces the classical graph-theoretical concepts by topological analogues. For example, instead of graph cycles one uses topological circles. This approach has been very fruitful, not only extending theorems from the finite to the infinite case (see e.g.\ \cite{LocFinMacLane,DiestelBook05,fleisch}), but also having further applications 
\cite{AgCurrents} and opening new directions \cite{basis,FundGp,DiSpHom,Hom2,ltop}. See \cite{RDsBanffSurvey} for a survey on this project. 

This paper is motivated by an attempt to generalise some of these graph-theoretical facts to continuous objects. And indeed, our main result is a generalisation of one of the most basic tools in the aforementioned project of Diestel \etal, \Tr{Zer} below, from graphs to arbitrary continua. In order to achieve this generalisation we introduce a new homology that generalises the cycle space of graphs to arbitrary metric spaces. We use this homology to conjecture a characterisation of the continua embeddable in the plane.

\subsection{Background and motivation}
The cycle space \ccg\ of a finite graph \g coincides with its first, simplicial or singular, homology group. As an example of the usefulness of this concept in graph theory, let me mention the following classical theorem of MacLane, providing an algebraic characterisation of the graphs embeddable in the plane.

\begin{theorem}[MacLane \ThmMacLaneC, \cite{DiestelBook05}] \label{mcl}
A finite graph \g is planar if and only if its cycle space \ccg\ has a 2-basis.
\end{theorem}
Here, a 2-basis is a set $B$ generating \ccg\ such that no edge of \g is used by more than two elements of $B$. See \cite{DiestelBook05} for more.

If the graph is infinite though, then \Tr{mcl} does not hold any more if \ccg\ is still taken to be the first simplicial or singular homology group \cite{LocFinMacLane}. However, \DK\ \cite{cyclesI,cyclesII} introduced a new homology for infinite graphs, called the \defi{topological cycle space \ccg}, which allows a verbatim generalisation of \Tr{mcl}:

\begin{theorem}[Bruhn \& Stein \cite{LocFinMacLane}] \label{mclinf}
A locally finite graph \g is planar if and only if its topological cycle space \ccg\ has a 2-basis.
\end{theorem}

The topological cycle space allows for such generalisations of all the fundamental facts about the cycle space of a finite graph. It is defined as a vector space, over $\Z_2$, consisting of sets of edges of the graph. Namely, it contains those edge-sets of \g that form topological circles in the end-compactification $|G|$ of \G, as well as the sums of these edge-sets, where we allow sums of infinitely many summands as long as they are well defined. An important innovation in the approach of \DK\ is that even if one is interested in the graph $G$ only, it is helpful to  consider the larger space $|G|$ that also contains the \defi{ends} of $G$.
The interested reader can find more  details and results about the  topological cycle space in \cite[Chapter~8.5]{DiestelBook05} or \cite{RDsBanffSurvey}; these details are however not necessary for understanding the current paper. The topological cycle space \ccg\ of \g is larger than the first simplicial homology group of \G, since the latter does not have any element comprising infinitely many edges. 

It is far less obvious, but true \cite{DiSpHom}, that \ccg\ is on the other hand smaller than the first singular homology group of \fcg. Consider for example the graph \g of \fig{kringel.COL}, which is a one-way infinite `ladder'. The end-compactification $|G|$ of \g is in this case its one-point compactification (graphs are considered as 1-complexes throughout the paper). Thus there is a loop \sig\ in $|G|$, depicted in \fig{kringel.COL}, starting at the top-left vertex $v$, winding around each of the infinitely many 4-gonal faces of \G, reaching the point at infinity, then returning to $v$, and finally winding around the whole graph once in the clockwise direction without using any of the perpendicular edges. It turns out \cite{DiSpHom} that \sig\ does not belong to the trivial element of $H_1(\fcg)$, but it does correspond to the trivial element of \ccg: it traverses each edge the same number of times in each direction; thus, seen as an element of \ccg, it is the empty set of edges. A similar example can be obtained in the Hawaiian earing by contracting a spanning tree of \G\ to a point. This pathological behaviour of \sig\ is due to the fact that although it winds around any hole the same number of times in each direction, it does so in such a complicated order that one cannot `disentangle' it by adding only finitely many boundaries of 2-simplices. To put it in a different way, the homology class of \sig\ is a product of infinitely many commutators.

\epsfxsize=0.6\hsize
\showFig{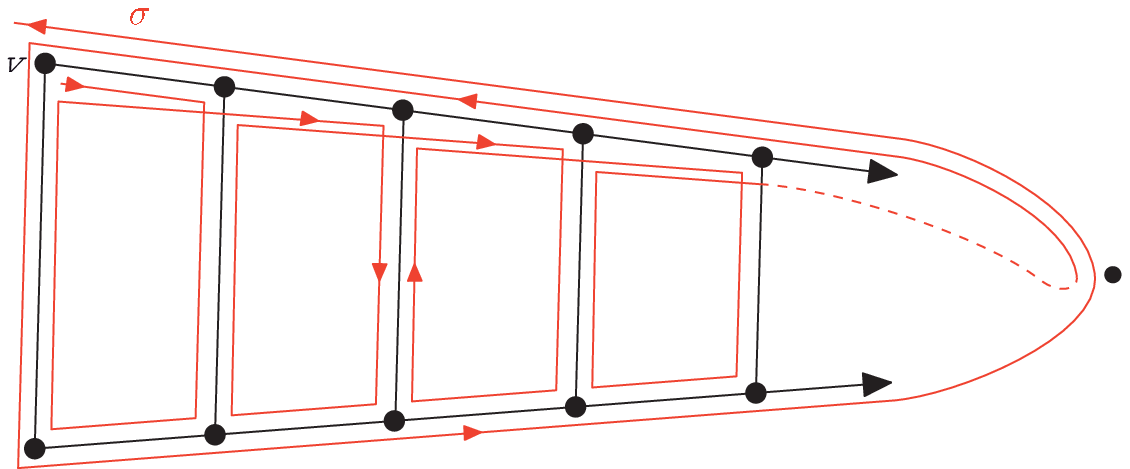}{A loop \sig\ that is not null-homologous although we would like it to be.}

This example shows that \ccg\ is indeed smaller than the first singular homology group of \fcg\ as claimed. However, this discrepancy between \ccg\ and $H_1(\fcg)$ should by no means be considered as a shortcoming of \ccg; for example, it is important for the truth of \Tr{mclinf}: the set of edge-sets of the 4-gonal faces of  \fig{kringel.COL} form a 2-basis, but it cannot represent a loop like \sig. It turns out, and is not hard to check, that \ccg\ is canonically isomorphic to the first \v{C}ech homology group of \fcg; see \cite{DiSpHom} for details.

We would like to generalise graph-theoretical theorems like \Tr{mcl} to continuous spaces. The main aim of this paper is to achieve such a generalisation for the following fact, which has been a cornerstone in the aforementioned project of Diestel \etal

\begin{theorem}[\DK\ \cite{cyclesII}]
  \label{Zer}
  Let $G$ be a  locally finite graph. Then every element of $\CCC(G)$ is a disjoint union of edge-sets of circles in \fcg.
\end{theorem}

\Tr{Zer} has found several applications in the study of \ccg\ \cite{LocFinMacLane,cyclesI,geo} and elsewhere \cite{fleisch}, and at least four proofs have been published; see \cite{hotchpotch} for an exposition.

Now in order to be able to generalise theorems like \Tr{Zer} or \Tr{mcl} to continuous spaces, we have to overcome two major difficulties: firstly, reformulate the assertions to rid them of any concepts, e.g.\ edges, that only make sense for graphs, and secondly, choose the right homology theory. 

To see how the first difficulty can be overcome, suppose that the graph $G$  in \Tr{Zer} is finite. We could then reformulate the assertion as follows: 
\labtequ{minl}{Every element of $\CCC(G)$ has a representative of minimal length.}
Here, a \defi{representative} is a formal sum of edge-sets of cycles. Indeed, this formulation is equivalent to that of \Tr{Zer} if \g is finite: a representative of minimal length cannot have two summands $C_1,C_2$ containing the same edge $e$, for then we could delete $e$, and any other common edges, from both $C_1,C_2$ and combine the remaining paths into a new cycle or new set of cycles whose total length is smaller, since we saved some length by removing $e$.

Formulation \eqref{minl} has the advantage that it makes sense for objects other than graphs if one replaces \ccg\ by some suitable homology group. The question now is, which homology should one use to extend this assertion beyond graphs. For example, singular homology will not do because of the example of \fig{kringel.COL}: the loop \sig\ has finite length if we metrize that space using the Euclidean metric, but there are loops homologous to \sig\ with arbitrarily small length, namely, those obtained by translating \sig\ to the right by one or more squares. Singular homology can fail to satisfy \eqref{minl} even if it is finitely generated, see \Er{excamera}.

\subsection{A new homology} \label{Snewh}

In view of the above discussion it is clear that in order to make assertion  \eqref{minl} true in general we need a homology group that excludes some `redundant' elements of singular homology. In fact such an approach is often followed when dealing with `wild' spaces, e.g.\ spaces with an uncountably generated fundamental group: in these cases many elements of the homotopy or homology groups do not capture some `hole' of the space but rather represent a complicated way to wind around infinitely many holes, and one wants to omit these elements in order to obtain a smaller group that still reflects the structure of the space; see \cite{Hom2,Dydak, EdaHom} for some examples. In many cases the better-known shape groups \cite{shape} also provide such simplifications of the corresponding homotopy or homology groups.


Earlier constructions of homology groups are not well-suited for our purposes as they either obviously fail to satisfy \eqref{minl} or it is not clear how to assign lengths to their representatives. In this paper, we will introduce a homology group \homp\ that comes with a natural notion of length, has the topological cycle space as a special case (\Sr{Sappg}) and, more importantly, makes assertion  \eqref{minl} true for all compact metric spaces.

We define \homp\ as a quotient of the first singular homology group $H_1$. For example, we would like to identify the class of \sig\ in the example of \fig{kringel.COL} with the trivial class. In order to decide which classes should be identified, we introduce a natural distance function on $H_1$, and identify any two elements if their distance is zero. This distance function is defined as follows. Intuitively, if two 1-cycles are not homologous, then there are some `holes' in our space that witness this fact, and we assign a distance to the corresponding pair of classes of $H_1$ reflecting the `size' of these holes. More precisely, the distance between two classes $c,d\in H_1$ is defined to be the minimal total area of a ---possibly infinite--- set of metric discs and cylinders that we could glue to our space $X$ to make $c$ and $d$ homologous. These metric discs and cylinders must bear a metric such that this glueing does not affect the metric of $X$. See \Sr{Sdeflhom} for the formal definitions. In \Sr{secEx} we display some examples that justify this definition by showing that modifying it would make assertion \eqref{minl} false.

An important feature of this distance function is that an infinite commutator product as the one of \fig{kringel.COL} can have distance zero to the trivial element. For example, patching all but finitely many of the 4-gonal faces in \fig{kringel.COL} by adding the missing trapeze would render \sig\ null-homologous, and this can be accomplished by adding arbitrarily little area if we skip a lot of the  4-gonal faces.

The aforementioned distance function gives rise to a metric on \homp\ after the identifications have taken place, which turns \homp\ into a metrizable topological group. We will also consider the completion \lhom\ of \homp, which will have the effect of strengthening our main result.


\subsection{Main result} \label{IntMain}
We can now state our main result. 

\begin{theorem}\label{mainI}
For every compact metric space $X$ and $C \in \hompx$, there is a \rep\ $\seq{z}$ of $C$ whose length is at most the infimum of the lengths of all representatives of $C$.
\end{theorem} 

Here, a \rep\ can intuitively be thought of as a sum of infinitely many 1-cycles $z_i$. Formally, a \rep\ of $C$ is defined as a sequence $\seq{z}$ whose initial subsequences give rise to a sequence $({\sum_{j\leq i} z_j})_{i\in\N}$ of 1-cycles the homology classes of which converge to $C$ \wrt\ the metric of \homp; see \Sr{Sdeflhom} for details. The  \defi{length} of a \rep\ is the sum of the lengths of the simplices in $z_i$, the latter lengths being defined in the standard way (see \Sr{gendefs}).

For example, consider the subspace $X$ of the real plane depicted in \fig{figkammI}. Let \sig\ be a closed 1-simplex $\sig: \unin \to X$ that traverses each of the infinitely many circles in this space precisely once and has finite length. Let $\beta\in H_1(X)$ denote the homology class of the 1-cycle $1\sig$. Note that for every representative of $\beta$ there is a further representative of smaller length, obtained by avoiding to traverse some of the perpendicular segments. Thus no representative achieves a minimum length. Still, \Tr{mainI} yields a \rep\ $\seq{z}$ of minimum length: let for example each $z_i$ be a closed simplex winding around the $i$th circle once in a straight manner.



\showFig{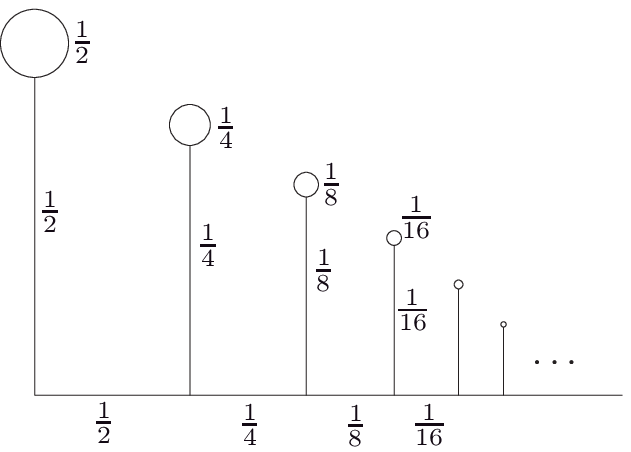}{A compact subspace of the real plane. The numbers denote the lengths of the corresponding segments.}


\medskip
\Tr{mainI} implies \Tr{Zer}. This can be shown by a similar argument as the one we used for the equivalence of the latter and \eqref{minl} for finite \G, except that if \g is infinite we assign lengths to its edges to make their total length summable; see \Sr{Sappg} for details. In fact, we obtain a strengthened version of \Tr{Zer}. Furthermore, with \Tr{mainI} we generalise \Tr{Zer} to \nlfg s, achieving a goal set by the author in \cite[Section 5]{ltop}; see \Sr{Sappg} for more.


For the proof of \Tr{mainI} we obtain an intermediate result which might be of independent interest. This result states that if $(H,+)$ is an abelian metrizable topological group, and a function $\ell: H \to \R^+$ is given satisfying certain natural properties that derive their intuition from the behaviour of lengths in geometry, then every element $h$ of $H$ can be `decomposed' as a sum $h = \sum h_i$ so that $\ell(h) = \sum \ell(h_i)$ and no $h_i$ can be decomposed further. See \Sr{Smezzo} for details.

\subsection{Further problems and remarks}

In this section we discuss some related conjectures for which strong evidence is available.

With \Tr{mainI} we extended a basic graph-theoretical tool to arbitrary compact metric spaces. It remains to try to exploit this in order to also extend results whose proofs are based on or related to this tool. A conjecture of this kind is offered in \cite[Conjecture 6.1]{ltop}. 
A further example is the following conjecture, which seeks an algebraic characterisation of the Peano continua embeddable in the plane, similar to that of \Tr{mcl}. 

\begin{conjecture}\label{conjMcL}
Let $X$ be a compact, locally connected, metrizable space that is locally embeddable in $S^2$. Then $X$ is embeddable in $S^2$ \iff\ \ti\ a simple set $S$ of circles in $X$ and  a metric $d$ inducing the topology of $X$ so that the set $U:= \{ \ec{\chi} \in \lhom (X,d) \mid \chi \in S\}$ {spans} $\lhom (X,d)$.
\end{conjecture}

See \cite[Conjecture 6.2]{ltop} for more on this conjecture. For example, $X$ here could be the Sierpinski triangle, in which case we could choose $S$ to be the set of its triangular face boundaries, corroborating the conjecture.

A further question motivated by our main result is whether something similar holds for higher dimensions. It is straightforward to see how to generalise the definition of \homp: instead of topological discs and cylinders one has to use their higher dimensional analogues. 
Our proof cannot prove this, but many of our intermediate steps still work. 

\begin{problem}
 Generalise \Tr{mainI} to higher dimensions.
\end{problem}
See \Sr{Shidi} for more on this problem.

Although we can generalise our homology group \homp\ or \lhom\ to higher dimensions, we do not obtain a homology theory in the sense of Eilenberg and Steenrod \cite{EilSteen, Hatcher}, since \hompx\ depends not only on the topology of $X$ but also on its metric. For the purposes of the current paper this is rather an advantage of \homp: since \Tr{mainI} holds for any choice of a compatible metric, we can affect the outcome of the application of the theorem by varying the metric. Still, it would be interesting to obtain a similar homology theory that does satisfy the axioms of Eilenberg and Steenrod by eliminating the dependence on the metric. Similarly, one could for example try to prove the following:

\begin{conjecture}
 Every Peano continuum $X$ has a metric compatible with its topology such that the corresponding \lhom\ coincides with the first \v{C}ech homology group of $X$.
\end{conjecture}

\Tr{lhiscc} below implies that this is true when $X$ is the end-compactification of a \lfg. The condition that $X$ be a Peano continuum is imposed because in a space that is not locally connected \v{C}ech homology may contain elements not represented by singular homology.

\numberwithin{equation}{section}

\section{General definitions and basic facts} \label{gendefs}

In this section we recall the standard definitions and facts that we will use later.
Most of this is very well-known but it is included for the convenience of the reader. For other standard terms used in the paper but not found in this section we refer to the textbooks \cite{armstrong} for topology, \cite{Hatcher} for algebraic topology and \cite{DiestelBook05} for graph theory.

For every metric space $M$, it is possible to construct a complete metric space $M'$, called the \defi{completion} of $M$, which contains $M$ as a dense subspace. The completion  $M'$ of $M$ has the following universal property \cite{MiVoFun}:

\begin{equation} 
\label{univ} 
\begin{minipage}[c]{0.85\textwidth} 
If $N$ is a complete metric space and $f: M \to N$ is a uniformly continuous function, then there exists a unique uniformly continuous function $f': M' \to N$ which extends $f$. The space $M'$ is determined up to isometry by this property (and the fact that it is complete).
\end{minipage}\ignorespacesafterend 
\end{equation}

Next, we recall the definition of the length of a topological path $\sigma: [a,b]\to M$ in a metric space $(M,d)$. For a finite sequence $S=s_1, s_2, \ldots, s_k$ of points in $[a,b]$, let $\ell(S):= \sum_{1\leq i< k} d(\sigma(s_i),\sigma(s_{i+1}))$, and define the \defi{length} of $\sigma$ to be $\ell(\sigma):=\sup_S \ell(S)$, where the supremum is taken over 
all finite sequences $S=s_1, s_2, \ldots, s_k$ with $a=s_1<s_2< \ldots <s_k=b$. This definitions coincides with that of the 1-dimensional Hausdorff measure of $\sig([a,b])$ when \sig\ is injective.

The $n$-dimensional Hausdorff measure of $M$ is defined by 
$$HM^n(M):= \lim_{\del \to 0} \inf \lbrace \sum_i  diam(U_i)^n \mid \bigcup_i U_i= M, diam(U_i) < \del \rbrace, $$
the infimum taken  over all countable covers \seq{U} of $M$ by sets $U_i$ of diameter less than \del.

We will also need the following.

\LemHeineCa{HeCa}

\section{Definitions and basic facts: \lhom, \rep s, and length; statement of main result} \label{Sdeflhom}


Let $X$ be any topological space, fixed throughout the paper, and consider its first singular homology group $\shom= \shom(X;\Gamma)$ over a group $\Gamma$. Our results are stated and proved for $\Gamma$ being any of the groups $\Z, \R$, or $\Z/n\Z$ for some $\nin$. We restrict ourselves to those groups because we want to make use of the absolute value $|a|$ of an element $a$ of $\Gamma$.

As mentioned in \Sr{Snewh} we want to put a distance function on $H_1$ and identify any two elements if their distance is zero. This distance between two classes $b,c$ measures the total area of the `holes' that we have to `patch' to make $b$ equivalent to $c$, in a sense that we will soon make precise. Intuitively, we are going to glue some spaces of a special form to $X$ in order to make $b$ equivalent to $c$, and measure the area of those spaces. Another way of saying that `we glue some spaces to $X$' is to say that `we embed $X$ into a larger space', and I found it more convenient to adhere to the second alternative. This motivates the following definition. 

\begin{definition} \label{darex}
An \defi{\arex} $(X',\iota)$ of $X$ is a metric space  $X'$ in which $X$ is embedded by an isometry $\iota: X \to X'$ such that each component of $X' \sm \iota(X)$ is either a \defi{metric disc} or a \defi{metric cylinder}, i.e.\ a metric space homeomorphic to either $\{x\in \R^2 \mid |x|<1\}$ or $\{x\in \R^2 \mid 1<|x|<2\}$. The \defi{\evol} of this \arex\ is the sum of the \vol\s\ of the components of $X' \sm \iota(X)$. 
\end{definition}
The \defi{\vol} of such a component can be defined as its 2-dimensional Hausdorff measure. However, in the \arex s that we will actually use, each such component is either a domain of $\R^3$ or a finite union of such domains. Thus the reader may choose,  for simplicity, to append to the definition of an \arex\ that each component of $X' \sm \iota(X)$ is a domain of  $\R^n$, and use his favourite definition of {\em area}; our proofs still apply without modification. I chose the above more generic definition because it leads to a stronger main result without complicating the proofs.

The effect of a metric disc in an \arex\ is to make a loop bounding it null-homologous. Similarly, the effect of a metric cylinder is to make two loops homologous to each other. Note that the latter could also be achieved by using two discs to make each of the two loops null-homologous. Thus one could wonder if we really need to allow for metric cylinders in \Dr{darex}. \Er{excamera} below shows however that we do need these metric cylinders in order to make our main result true.

\medskip
We now define a pseudo-metric $d_1$ on the singular homology group $H_1(X)$ of $X$. Given two elements $[\phi],[\chi]$ of $H_1(X)$, where $\phi$ and $\chi$ are $n$-chains, let $d_1([\phi],[\chi])$ be the infimum of the \evol\s\ of all \arex s $X'$ of $X$ \st\ $\phi$ and $\chi$ belong to the same element of $H_1(X')$.

It follows easily by the definitions that 

\labtequ{trinq}{$d_1$ satisfies the triangle inequality.}

However, $d_1$ is not yet a metric, since there may exist $c\neq f\in H_1$ with $d(c,f)=0$: for example, the homology class $c$ of the loop of \fig{kringel.COL} satisfies $d(c,\zero)=0$ although $c\neq \zero$ as proved in \cite{DiSpHom}. Still, declaring $c,f\in H_1$ to be equivalent if $d(c,f)=0$ and taking the quotient with respect to this equivalence relation we obtain the group $\homp=\hompx$; the group operation on $\homp$ can be naturally defined for every $c,d\in \homp$ by choosing representatives $\alp\in c$ and  $\beta\in d$ and letting $c+d:= \llbracket \alp+\beta \rrbracket$ be the class in \homp\ containing the element $\alp+\beta$ of $H_1$. To see that this sum is well defined, i.e.\ does not depend on the choice of $\alp$ and $\beta$, note that the union of two extensions of $X$ of \evol\ at most \eps\ each is an extension of $X$ of \evol\ at most $2\eps$.

We will use the notation $\llbracket \gamma \rrbracket$, where $\gamma$ is either an element of $H_1(X)$ or a 1-cycle, to denote the equivalence class in \hompx\ containing $\gamma$ or $[\gamma]$ respectively, where $[\chi]$ always denotes the element of $H_1(X)$ containing the 1-cycle $\chi$.

Now $d_1$ induces a distance function on $\homp$, which we will, with a slight abuse, still denote by $d_1$: for any $\ec{\phi},\ec{\chi}\in \homp$ let 
$d_1(\ec{\phi},\ec{\chi}) := d_1([\phi],[\chi])$;
it is an easy consequence of \eqref{trinq} that this is well defined, and that $d_1$ is now a metric on $\homp$.

\begin{definition} \label{deflhom}
We now define a new homology group $\lhom=\lhomx$ of $X$ to be the completion of \homp\ \wrt\ the metric $d_1$. The operation of $\lhom$ is defined, \fe\ $C,D\in \lhom$, by $C+D:= \lim_i (c_i + d_i)$ where \seq{c}\ is a \Cs\ in $C$ and \seq{d}\ is a \Cs\ in $D$.
\end{definition}

If $X$ is compact and well-behaved then it might be the case that $\lhom$ is complete, which means that it coincides with  \homp. However, this is not always the case; see \Er{homneq}. If $C\in \homp$ then we will sometimes, with a slight abuse of notation, still use the symbol $C$ to denote the element of \lhom\ corresponding to $C$, that is, the equivalence class of the constant sequence $(C)$.

\lhom\ is by definition a metrizable abelian topological group. If defined over \R\ then it can also be viewed as a Banach space.

The operation $C+D$ in \Dr{deflhom} is well defined since, by \eqref{trinq}, $c_i + d_i$ is a \Cs\ too and it does not depend on the choice of $c_i$ and $d_i$.

The following observation, which is easy to prove, can be used to obtain an alternative definition for the addition operation $C+D$, where one first adds 1-cycles and then considers their homology classes rather than the other way round. Here $(\phi_i)$ and $(\chi_i)$ are sequences of 1-cycles.

\begin{lemma}\label{sum}
Let $(\ec{\phi_i})$ and $(\ec{\chi_i})$ be \Cs\s\ in \homp. Then $\lim \ec{\phi_i +\chi_i} = \lim (\ec{\phi_i} + \ec{\chi_i})$.
\end{lemma} 
\noproof

\comment{
	The following basic fact follows easily from the definitions and \eqref{trinq}.

	\begin{lemma}\label{dtriang}
Let $X$ be a metric space and let $C,D\in \lhom(X)$. Then $d_1( (C+D) , \zero) \leq d_1( C, \zero) +d_1(D,\zero)$ (and thus $d_1( (C-D) , \zero) \geq d_1( C, \zero) -d_1(D,\zero)$).
\end{lemma} 
	\noproof
}

Before we go on to prove our main result about \lhom\ we should pause to think whether we just identified all elements of $H_1$ to the zero element to obtain a trivial \lhom, which would make our main result void. In fact, this can happen in certain pathological spaces, but we will show that, for example, $\lhom(S^1)$ is not trivial, and this can be applied to show the non-triviality of \lhom\ for many other spaces. See \Sr{secEx} for more.

\medskip

A \defi{\rep} of $C\in \lhom$ is an infinite sequence $\seq{z}$ of $1$-cycles $z_i\in Z_1$  such that the sequence $(\ec{\sum_{j\leq i} z_j})_{i\in\N}$ is a \Cs\ in $C$. One can think of a \rep\ as an ``1-cycle'' comprising the infinitely many 1-simplices $z_i$.  Later on (\Sr{Infsums}) we will rigorously define infinite sums of elements of \lhom, and it turns out that $C = \sum_i \ec{z_i}$ for every representative $(z_i)$ of $C$. 

For example, in \fig{figkamm}, we can build a \rep\ of the class of the loop described there by letting $z_i$ be an 1-simplex going around the $i$th circle.

One of the central concepts in our main \Tr{lCycD} is the \defi{length} of an element of $\lhom$. To define this we first need to define the length of a simplex, a 1-chain, and an element of \shom. With a considerable abuse of notation, we will denote the length of any of those objects by $\ell()$.

Since a simplex $\chi$ is by definition a topological path, we can use the standard definition of its length $\ell(\chi)$ as in \Sr{gendefs}. We can then define the length of a 1-chain $z= \sum_i a_i \chi_i$ by 
$\ell(z):= \sum_i |a_i| \ell(\chi_i)$, and consequently the length $\ell(\beta)$ of an elements $\beta$ of \shom\ by $\ell(\beta):= \inf_{z\in \beta} \ell(z)$. Finally, for $C \in \lhom$, we define $\ell(C):= \inf_{\beta_i} \lim_i \ell(\beta_i)$ where the infimum ranges over all sequences \seq{\beta}\ with $\beta_i\in \shom$ \st\ $(\llbracket\beta_i \rrbracket)_{\iin}$ is a \Cs\ in $C$ and $\lim_i \ell(\beta_i)\in \R^+ \cup \{\infty\}$ exists. 

We can now state our main result, \Tr{mainI}, in a stronger and more precise  form. Recall that an \defi{1-simplex} is a continuous function $\sig: [a,b] \to X$. If $\sig(a)=\sig(b)$ then \sig\ is called a \defi{closed simplex}, and if moreover \sig\ is injective on $[a,b)$ then it \sig\ is called a \defi{\clex} (note the similarity to a \defi{circle}, i.e. a homeomorph of $S^1$).

\begin{theorem}\label{lCycD}
For every compact metric space $X$ and $C \in \lhomx$, there is a \rep\ $\seq{z}$ of $C$ with $\sum_i \ell(z_i)=\ell(C)$. 

In particular, for every other \rep\ $\seq{w}$ of $C$ we have $\sum_i \ell(z_i)\leq \sum_i \ell(w_i)$.

Moreover, if $\ell(C)<\infty$ then \seq{z} can be chosen so that each $z_i$ is 
 a \clex. 
\end{theorem} 

As a consequence we can now simplify the definition of $\ell(C)$ once we have proved \Tr{lCycD}: the following assertion yields an equivalent definition.

\begin{corollary}
For every compact metric space $X$ and $C \in \lhomx$ we have $\ell(C) = \inf \sum_i \ell(z_i)$, the infimum ranging over all \rep s \seq{z}\ of $C$.
\end{corollary}
\begin{proof} 	
\Tr{lCycD} immediately yields $\ell(C) \geq \inf \sum_i \ell(z_i)$. The reverse inequality follows from the definition of $\ell(C)$: given any \rep\ \seq{z}\ of $C$ we can let $\beta_i:= [ \sum_{j\leq i} z_j]$, and as $(\llbracket\beta_i \rrbracket)_{\iin}$ is a \Cs\ in $C$ by the choice of $(z_i)$, we have $\ell(C) \leq \lim_i \ell(\beta_i) \leq \lim_i \ell(\sum_{j\leq i} z_j) = \sum_i \ell(z_i)$.
\end{proof}

\section{Isoperimetric properties of lengths} \label{Siso}

In this section we prove two basic facts relating length and area in metric spaces. The reader is encouraged to skip this section during the first reading of the paper and come back when it becomes relevant.

The following lemma yields a kind of isoperimetric property for arbitrary metric spaces: it shows that any ``hole'' can be ``filled in'' by an area proportional to the square of the perimeter of the hole.

\begin{lemma}\label{U}
There is a universal constant $U$ \st\ \fe\ metric space $(X,d_X)$ and every closed curve $\sigma:I \to X$ \ti\ an \arex\ $(X',\iota)$ of \evol\ at most $U \ell^2(\sig)$ in which \sig\ is \nuho. 

Moreover, $X'$ can be chosen so that $X' \sm \iota(X)$ is a \mts\ with diameter less than $\ell(\sig)$.
\end{lemma} 
\begin{proof} 	
Pick a (geometric) circle $D$ of length $\ell:=\ell(\sig)$ in $\R^3$, and a continuous mapping $f$ from $D$ to $I$ \st\ corresponding subpaths have equal lengths; that is, for any subarc $D'$ of $D$ we have $\ell(D')=\ell(\sigma\restr f(D'))$. Let $S$ be a closed hemisphere in $\R^3$ having $D$ as its equator, and give $S$ its path metric (i.e.\ the distance $d_S$ of two points in $S$ is defined to be the minimum length of an arc in $S$ between these two points). Now in order to obtain the desired \arex\ $X'$, we glue a copy of $S$ along the image of \sig\ in $X$ using $\sig \circ f$ as an identifying map. We still have to specify a metric $d'$ for $X'$. Note that by the choice of $d_S$, for every pair of points $x,y$ in the domain of $\sig \circ f$ we have 
\labequ{dxs}{d_S(x,y)\geq d_X(\sig \circ f(x), \sig \circ f(y)).}
This allows us to extend the metric $d_X$ of $X$ into a metric $d'$ of $X'$ as follows. Let $d'(z,w)= d_X(z,w)$ for every pair of points $z,w$ of $X$, including points that got identified with points of $S$. If $z\in S$ and $w\in X$, then let $d'(z,w)= \inf_y \{d_S(z,y) + d_X(\sig \circ f(y),w)\}$, the infimum taken over all points $y$ in $D$. Finally, if $z,w\in S$ let $$d'(z,w)= \min \{ d_S(z,w) , \inf_{y,y'} \{d_S(z,y) + d_X(\sig \circ f(y), \sig \circ f(y')) + d_S(y',w)\} \},$$ the infimum taken over all pairs of points $y,y'$ in $D$, even if $y=y'$. It is an easy exercise to check, using \eqref{dxs}, that $d'$ is indeed a metric, and that $(X',id)$ is an \arex\ of $X$ of \evol\ at most the area of $S$. 

The Euclidean area of $S$ is $2\pi R^2$ for $R:=\ell/2\pi$. Since we consider the path metric $d_S$ on $S$, distances are greater by a factor of up to $\pi$ compared to the Euclidean metric, thus $area(S)\leq 2\pi R^2 \pi^2=\pi \ell^2/2$ and we can take $U=\pi/2$. Moreover, the diameter of $S$ is by construction $\ell(\sig)/2$. This completes the proof.

\medskip
A point that might require some clarification is that we are not assuming that the closed curve \sig\ is injective in its interior. If it is not, then the closure of $X'\sm X$ is not necessarily a closed disc, but still $X'\sm X$ itself is an (open) metric disc as the reader can check, and so $X'$ is indeed an \arex.
\end{proof}

The following observation about the real numbers is an easy exercise.
\labtequ{squares}{For every $\ell, \eps \in \R^+$ \ti\ an $r\in \R$ \st\ if $a_1, a_2, \ldots, a_k$ are positive real numbers with $a_i<r$ \fe\ $i$ and  $\sum a_i = \ell$, then $\sum a_i^2 < \eps$.}


Our previous lemma shows that a `hole' of small perimeter can be patched using relatively little area. Our next result performs a similar task: it shows that if two holes are bounded by curves that are `close' to each other, then the corresponding homology classes can be made equivalent using relatively little area. The following definition makes this concept of `closeness' precise; see  also \fig{cylinder}.

\begin{definition} \label{dclose}
Let $\sigma,\tau:I \to X$ be two closed curves in a metric space $X$. We will say that $\sig$ and $\tau$ are \defi{$\del$-close}, if $|\ell(\sig) - \ell(\tau)|<\del$ and moreover there are subdivisions $\sig^1, \sig^2, \ldots \sig^k$ and $\tau^1, \ldots \tau^k$ of $\sig$ and $\tau$ respectively that fulfill the following requirements \fe\ $i\in \{1,\ldots,k\}$: 
\begin{enumerate}
\item $\ell(\sig^i)< \del$ and $\ell(\tau^i) <\del$,
\item \label{cliii} $|\sum_{j\leq i} \ell(\sig^j) - \sum_{j\leq i} \ell(\tau^j)| < \del$, and
\item if $p,q$ are the vertices of $\sig^i$ and $p',q'$ are the vertices of $\tau^i$ then $d(p,p')<\del/k$ and $d(q,q')<\del/k$.
\end{enumerate}
\end{definition}

\showFig{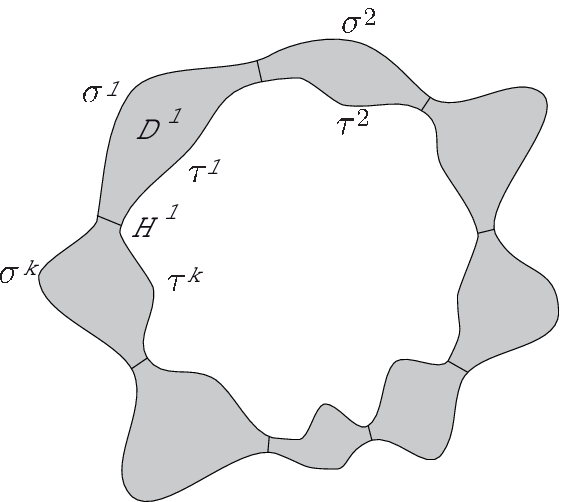}{Two curves that are close to each other can be made homotopic by a cylinder of relatively small area.}

Define the \defi{width} of a homotopy $h: \unin \times I \to X$ to be $$\sup \{ d( h((0,y)), h((x,y)) ) \mid x\in \unin, y\in I\}.$$

We say that a curve $\sigma:I \to X$ has \defi{\cosp} $c$ if for every subinterval $[a,b]$ of $I$ we have $\frac{\ell(\sig([a,b]))}{b-a} = c$. 

We can now state our next lemma.

\begin{lemma}\label{Uc}
For every $\eps,l \in \R_+$ with $\eps <l$ there is an $f(l,\eps)\in \R_+$ \st\ \fe\  metric space $X$, and every two closed curves $\sigma,\tau:I \to X$ in $X$ of length less than 
$l+\eps$  
that are $f(l,\eps)$-close, there is an extension $X'$ of $X$ of \evol\ less than \eps\ in which \ti\ a homotopy $h$ between $\sig$ and $\tau$. Moreover, if \sig\ and $\tau$ have \cosp\ then $X'$ and $h$ can be chosen so that the width of $h$ is less than $5f(l,\eps)$.
\end{lemma} 
\begin{proof} 	
Suppose the closed curves $\sigma,\tau$ of length less than $l+\eps$ are $\del$-close for some real number $\del$ much smaller than $l$, and let $\sig^1, \ldots \sig^k$ and $\tau^1, \ldots \tau^k$ be subdivisions of $\sig$ and $\tau$ respectively as in \Dr{dclose}. In order to construct the desired extension $X'$, start by adding to $X$, \fe\ $1\leq i \leq k$, an isometric copy $H^i$ of the real interval $[0, m^i]$, where $m^i$ is the distance in $X$ between the first vertex $p_i$ of $\sig^i$ and the first vertex $p'_i$ of $\tau^i$; then identify one endpoint of $H^i$ with $p_i$ and the other with $p'_i$, see \fig{cylinder}. After having done so, note that concatenating,  \fe\ $i$, the ``paths'' $\sig^i, H^{i+1}$, $\tau^i$ (inversed) and $H^i$ we can obtain a closed curve $c^i $ of length 
\labtequ{lci}{$\ell(c^i)< \del(2 + 2/k)$.}
As in the proof of \Lr{U}, we can \fe\ $i$ glue a disc $D^i$ of small diameter along the image of $c^i$ to obtain an extension of $X$ of \evol\ at most $U \ell^2(c^i)$ in which $c^i$ is contractible. Uniting all these extensions ---identifying points corresponding to the same point of $X$--- we obtain the extension  $X'$ of \evol\ $V\leq U\sum \ell^2(c^i)$. It is easy to see that $\sig$ is indeed homotopic to $\tau$ in $X'$. Note that 
$$\sum \ell(c^i) < \ell(\sig) + \ell(\tau) + 2\sum \ell(H^i) \leq 2 (l + \eps + \del)< 3l,$$ and as the length of each $c^i$ is bounded from above by~\eqref{lci}, it follows from~\eqref{squares} that choosing $\del$ small enough we can achieve $V<\eps$ as desired. Thus we can let $f(l,\eps):= \del$ for such a \del.

\medskip
To prove the second sentence of the assertion, suppose now that \sig\ and $\tau$ have \cosp, and define $X'$ as above. We are now going to construct the desired homotopy $h$. To begin with, let $h(0,x) = \sig(x)$ and $h(1,x)= \tau(x)$. Moreover, \fe\ $1\leq i \leq k$ let $h$ map the straight line segment $L^i$ in $\unin \times \unin$ joining the preimages of the first vertices of $\sig^i, \tau^i$ homeomorphically to $H^i$. Note that the segments $L^i$ do not intersect each other except perhaps at their endpoints. Then, extend $h$ continuously to the rest of  $\unin \times \unin$, mapping the area bounded by $L^i$ and $L^{i+1}$ to the disc $D^i$.

We claim that the width of $h$ is less than $5f(l,\eps)$. To see this, consider a point $p=(x,y) \in \unin \times \unin$, let $D$ be the disc containing $h(p)$ (or one of the discs containing $h(p)$ if we were unlucky and $h(p)$ lies in some $H^i$) and let $D'$ be a disc whose boundary contains $h(p')$ where $p':=(0,y)$. Requirement~\ref{cliii} of \Dr{dclose} provides a lower bound for the angles that the segments $L^i$ form with the segment $\{0\} \times \unin $; more precisely, requirement~\ref{cliii} and the fact that $\sig$ and $\tau$ have \cosp\ implies that there is a point $t=(0,t_1)\in \unin \times \unin$ \st\ $h(t)\in D$ and the length $\mu$ of the restriction of $h(0,x)=\sig(x)$ to the interval between $t$ and $p'$ is at most $2f(l,\eps)$. Indeed, let $L^m$ be the segment separating $p$ from $p'$ in $\unin \times \unin$, and let $t=(0,t_1)$ and $t'=(0,t'_1)$ be the endpoints of $L^m$ in $\{0\} \times \unin $ and $\{1\} \times \unin $ respectively. By the choice of $L^m$ we have $|t_1 - t'_1| \geq |t_1 - y |$. But recall that $t,t'$ are the preimages of first vertices of $\sig^m, \tau^m$ respectively, and since $\sig$ and $\tau$ have \cosp, and by~\ref{cliii} of \Dr{dclose} there holds  $|\sum_{j\leq m} \ell(\sig^j) - \sum_{j\leq m} \ell(\tau^j)| < f(l,\eps)$, it follows that the length $\mu$ of the restriction of $h(0,x)=\sig(x)$ to the interval between $t$ and $p'$ is at most $2f(l,\eps)$ as claimed.

By the definition of length, this implies $d(h(t), h(p'))\leq 2f(l,\eps)$. Moreover, as both $h(t), h(p)$ lie in $D$, and  $D$ was chosen so that its diameter is at most $3f(l,\eps)$, we have $d( h(p),h(p') )\leq 3f(l,\eps) + 2f(l,\eps)= 5f(l,\eps)$. Since $p$ was chosen arbitrarily, the last inequality proves that the width of $h$ is at most $5f(l,\eps)$.

\end{proof}

\section{Basic facts about lengths}

In this section we prove some basic facts about lengths of homology classes, as defined in \Sr{Sdeflhom}, which we will need later.

Our first task is to prove that  $\ell(C)$ is attained by some sequence $(\beta_i)$ \fe\ $C \in \lhom$:

\begin{observation} \label{obatt}
For every $C \in \lhom$ \ti\ a sequence \seq{\beta}\ with $\beta_i\in \shom$ \st\ $(\llbracket \beta_i \rrbracket)_{\iin}$ is a \Cs\ in $C$ and $\ell(C) = \lim_i \ell(\beta_i)$.
\end{observation}

Note that this observation follows immediately from our main result \Tr{lCycD}, by taking $\beta_i = [\sum_{j\leq i} z_j]$, but \Tr{lCycD} is much stronger. As we will use \Or{obatt} in the proof of \Tr{lCycD} we have to prove the former separately:

\begin{proof}[Proof of \Or{obatt}]
If $\ell(C)=\infty$ then the assertion is easily seen to be true, so suppose $\ell(C)<\infty$. Pick a non-constant sequence \seq{h}, $h_i\in \R$, that converges to $\ell(C)$ from above. For $j=1,2,\ldots$, let \seq{\beta^j}\ be a sequence of elements of $\shom$ \st\ $(\llbracket\beta^j_i\rrbracket)_{\iin}$ is a \Cs\ in $C$ and $\lim_i \ell(\beta^j_i) < h_j$; such a sequence exists by the definition of $\ell(C)$. Pick an index $k\in \N$ \st\ $\ell(\beta^j_k) < h_j$ and $d_1(C,\llbracket\beta^j_k\rrbracket)< 2^{-j}$, and let $\beta_j:= \beta^j_k$; such a $k$ exists by the choice of $(\beta^j_i)$. 

By construction, the sequence \seq{\beta}\ we just constructed is a \Cs\ in $C$ and satisfies $ \lim_i \ell(\beta_i) = \ell(C)$.
\end{proof}

Since for every $\beta\in\shom$ we can, by the definition of $\ell(\bet)$, find 1-cycles in $\bet$ with lengths arbitrarily close to $\ell(\bet)$, we obtain with the above observation

\begin{corollary} \label{obl}
For every $C \in \lhom$ \ti\ a sequence of 1-cycles \seq{\phi}\ \st\ $(\ec{\phi_i})_{\iin}$ is a \Cs\ in $C$ and $\ell(C) = \lim_i \ell(\phi_i)$.
\end{corollary}

With \Lr{U} \Cr{obl} easily yields

\begin{observation} \label{zerol}
If $C\neq \zero \in \lhom$ then $\ell(C)>0$.
\end{observation}

Next, we check that the lengths of elements of $H_1$ satisfy a triangle inequality:

\begin{lemma}\label{triang}
Let $X$ be a metric space and let $\phi,\chi$ be two 1-chains in $X$. Then $\ell([\phi+\chi])\leq \ell([{\phi}])+\ell([{\chi}])$ (and thus $\ell([\phi-\chi])\geq \ell([{\phi}])-\ell([{\chi}])$).
\end{lemma} 
\begin{proof}
It is a trivial fact that if $\phi',\chi'$ are 1-chains in $X$, then $\ell(\phi'+\chi') \leq \ell({\phi'})+\ell({\chi'})$. The assertion now easily follows from the definition of $\ell([\phi])$, since if $b_1,b_2 \in B'_1$ then $\ell(\phi+\chi+b_1+b_2) \leq \ell({\phi+b_1})+\ell({\chi+b_2})$.
\end{proof}


From this we easily obtain a triangle inequality for elements of \lhom\ too:

\begin{corollary}\label{ltriang}
Let $X$ be a metric space and let $C,D\in \lhom(X)$. Then $\ell(C+D)\leq \ell(C)+\ell(D)$ (and thus $\ell(C-D)\geq \ell(C)-\ell(D)$).
\end{corollary} 
\noproof

\section{Examples} \label{secEx}

In this section we show examples that explain some of our choices in the preceding definitions and statements.

We defined \homp\ as a quotient of $H_1$ by identifying pairs of elements with distance zero. This identification entails the danger of identifying all of $H_1$ with the trivial element, which would make \homp\ and \lhom\ trivial and our main result void. And indeed, in certain pathological spaces $X$, e.g.\ when each element of $H_1(X)$ can be represented as an infinite product of commutators, this could happen; an example of such a spece can be found in  \cite{BoZaInf}. However, the following basic example when $X= S^1$ shows that \homp\ and \lhom\ are not trivial when it should not be:

\begin{theorem}\label{TSi}
$\lhom(S^1) \cong H_1(S^1)$.
\end{theorem} 
\begin{proof} 	
Let \sig\ be a \clex\ in $S^1$. Then $H_1(S^1)$ is generated by the corresponding homology class $[\sig]$. Thus all we need to show is that $[\sig]$ is not identified with the trivial element $\zero$ of $H_1(S^1)$; in other words, that \ti\ a lower bound $M$ \st\ every \arex\ of $S^1$ in which \sig\ is null-homologous has \evol\ at least $M$. 

So let $(S',\iota)$ be an \arex\ of $S^1$ in which \sig\ is null-homologous. Thus there is a 2-chain $\cb$ in $S'$ whose boundary is \sig. From now on we assume for simplicity that the group of coefficients on which $H_1(S^1)$ is based is \Z; the interested reader will be able to adapt our arguments to other groups of coefficients.

We may assume \obda\ that $\cb$ consists of a single 2-simplex $\rho$ whose boundary is a subdivision of \sig\ into three subsimplices, for otherwise we can combine pairs of 2-simplices of \cb\ together to get a shorter 2-chain. All we need to show now is that the area $A(P)$ of the image $P$ of $\rho$ is bounded from below by some constant $M$ independent of $S'$. In fact, we will show that we can  choose $M=1$.
Recall that we defined the area of a metric space to be its 2-dimensional Hausdorff measure, although the reader could probably arrive to the same conclusions using any alternative concept of area he is keen on.
Note that this bound $M=1$ is best possible, since it equals the 2-dimensional Hausdorff measure of the unit disc.

To prove the claimed bound for $A(P)$, we subdivide $S^1$, and \sig, into four equal arcs $X_1,Y_1,X_2,Y_2$ of length $\pi/2$ each, traversed by \sig\ in that order. Note that 
\labtequ{dsq}{$d(x,X_1) \geq \sqrt{2}$ holds for every $x\in X_2$, and similarly for $Y_1,Y_2$.}
We will use this observation to prove that $A(P)$ is at least half the area of a square with side length $\sqrt{2}$.

For this, define the mapping $f: P \to [0,\sqrt{2}]^2$ by $x\mapsto ( \ceil{d(x,X_1)},\ceil{d(x,Y_1)} )$ where $\ceil{d}:= \max \{ d, \sqrt{2}\}$. We may assume that the domain of $\rho$ is also the square $[0,\sqrt{2}]^2$ rather than the standard 2-simplex. Now consider the function $f \circ \rho: [0,\sqrt{2}]^2 \to [0,\sqrt{2}]^2$, which is continuous since both $f$ and $\rho$ are. Note that the restriction of $f \circ \rho$ to the boundary of the square $[0,\sqrt{2}]^2$ is, by \eqref{dsq}, a homeomorphism from that boundary onto itself. From this we will infer that 
\labtequ{frho}{$f \circ \rho$ is onto.}
There are perhaps many ways to prove this basic fact, and the reader might have a favourite one depending on their background. Here we sketch a proof using homology: suppose, to the contrary, that some point $z\in [0,\sqrt{2}]^2$ is not in the image $I$ of $f \circ \rho$, and let $Q:= [0,\sqrt{2}]^2 \sm \sgl{z}$. Note that $Q$ is homotopy equivalent to $S^1$, and so $H_1(Q)$ is isomorphic to $H_1(S^1)$ \cite[Corollary 2.11]{Hatcher}. But $f \circ \rho$ is a 2-simplex of $Q$ proving that its boundary is null-homologous, and so $H_1(Q)$ is trivial by the remark preceding \eqref{frho}. This contradiction establishes \eqref{frho}. Note that $f$ must thus also be onto.

Now suppose that $A(P)< 1$, which means that \fe\ \del\ there is a countable cover $\seq{U}$ of $P$ with $diam(U_i)< \del$ and $\sum diam(U_i)^2<1$. Letting $V_i:= f(U_i)$ we obtain a cover $\seq{V}$ of $[0,\sqrt{2}]^2$ since $f$ is onto. Moreover, by the definition of $f$ and the triangle inequality we have $diam(V_i)\leq \sqrt{2} diam(U_i)$. Thus $\sum diam(V_i)^2 \leq 2 \sum diam(U_i)^2< 2$ by the above assumption. This means that the area of $[0,\sqrt{2}]^2$ is less than 2, a contradiction.

This completes the proof that any \arex\ of $S^1$ in which \sig\ is null-homologous has an \evol\ of at least 1, implying that $\lhom(S^1) \cong H_1(S^1)$.
\end{proof}

Using the same arguments one can generalise this to the following.

\begin{corollary}\label{c}
Let \g be a \lf\ 1-complex. Then $\lhom(G) \cong H_1(G)$.
\end{corollary} 
\noproof

The following important example shows that \Tr{lCycD} would become false if we banned metric cylinders from the definition of an \arex. Moreover, it shows that \Tr{lCycD} fails if we replace \lhom\ by the first singular homology group $H_1(X)$ even if $H_1(X)$ is finitely generated. This example could also contribute to a better understanding of \Sr{Sprim}, the geometric part of the proof of our main result.

\example{excamera}{
\normalfont
We will define our space $X$ as a subspace of $\R^3$ with the Euclidean metric. It is similar to a well-known construction of \cite{archipelago} called the \defi{harmonic archipelago}.
The shape of $X$ is reminiscent of the shape of an old-fashioned folding camera: for every even $i\in \N$ let $D_i$ be the circle $\{(x,y,z)\in \R^3 \mid y= 2^{-i}, x^2 + z^2 = 1 \}$ and \fe\ odd $i\in \N$ let $D_i$ be the circle $\{(x,y,z)\in \R^3 \mid y= 2^{-i}, x^2 + z^2 = 1/2 + 2^{-i} \}$. Moreover, \fe\ \iin\ let $X_i$ be the closed cylinder in $\R^3$ with boundary $D_i \cup D_{i+1}$ that has minimum area among all such cylinders. Let $X$ be the closure of $\bigcup X_i$ in $\R^3$, that is, $X$ is the union of $\bigcup X_i$ and the cylinder $\{(x,y,z)\in \R^3 \mid y= 0, 1/2 \leq x^2 + z^2 \leq 1 \}$. 

For every $\iin$ let $\sig_i$ be a \clex\ that travels once around $D_i$. Note that $\sig_i$ is homotopic to $\sig_j$ \fe\ $i,j$. However, no $\sig_i$ is homologous to a \clex\ $\tau$ that travels once around the circle $D:=\{(x,y,z)\in \R^3 \mid y= 0, x^2 + z^2 = 1/2 \}\subset X$, because no 2-simplex can meet infinitely many $X_i$. Moreover, the two homology classes corresponding to $\tau$ and the $\sig_i$ cannot be made equivalent by glueing discs of arbitrarily small area to $X$ without distorting its metric. Thus, if we modified the definition of \arex\ to only allow discs as components of $X' \sm \iota(X)$, then \Tr{lCycD} would fail for $C:= \ec{1\sig_1}$, as $C$ has representatives with length arbitrarily close to $\pi=\ell(\tau)$, namely, the $\sig_i$, but no representative of length $\pi$ or less.

This example also shows that we cannot replace \lhomx\ by  $H_1(X)$ (and `\rep' by `representative') in the assertion of \Tr{lCycD} even if $H_1(X)$ is finitely generated. Indeed,  $H_1(X)$ is generated by 2 elements here, namely $[\sig_1]$ and $[\tau]$, and $[\sig_1]$ has no representative of minimum length.
\noproof

}

If $X$ is compact then in many cases we do not gain anything when we take the completion \lhomx\ of \hompx. For example, if $X$ is the space of \fig{figkamm} then $\hompx$ is already complete as the interested reader can check. There are however compact examples $X$ where $\hompx$ is not complete:

\example{homneq}{
\normalfont
Let $X$ be a metric space obtained as follows. Start with a topologist's sine curve $S$, pick a countably infinite `cofinal' sequence \seq{u}\  of points of $S$, and attach a circle of length $2^{-i}$ at each point $u_i$. To see that \hompx\ is not complete, let $\sig_i$ be a \clex\ corresponding to the circle attached at $u_i$, and note that $(\llbracket\sig_i\rrbracket)_\iin$ is a \Cs\ that has positive distance from each element $c$ of $H_1(X)$. Indeed, any such $c$ must miss some circle, and \Tr{TSi} yields a lower bound for that distance. 

\noproof
}

\comment{
The following is false in general, since an element of $\homp$ can contain two classes one of which has infinite length and the other not (Philipp and Reinhard's example)

\begin{observation}
If \seq{\beta}\ and \seq{\alpha}\ are two \Cs s in $c$, then $lim_i \ell(\beta_i) =lim_i \ell(\alpha_i)$.
\end{observation}
\begin{proof}
These limits exist since $d_1$ is a metric on $\homp$. Thus it suffices to show that $lim_i (\ell(\phi_i) - \ell(\chi_i))=0$. But since $d-1$ is a metric, we have $\ell(\phi_i) - \ell(\chi_i) \leq \ell(\phi_i-\chi_i)=d(\phi_i,\chi_i)$, and as the sequences \seq{\phi}\ and \seq{\chi}\ are equivalent, this amount does converge to $0$.
\end{proof}
}

For $C' \in \homp$ the element $C$ of \lhom\ corresponding to $C'$  satisfies $\ell(C)\leq \inf \{\ell(\beta) \mid \beta \in C'\}$ by the definitions. The aim of our next example is to show that this inequality can be proper. This means that (the first sentence of) the assertion of \Tr{lCycD}, applied to a $C\in \homp$, is in fact stronger than that of \Tr{mainI}.

\example{kamm}{
\normalfont
Consider the compact space $X\subseteq \R^2$ depicted in \fig{figkamm}. It is easy to construct a closed 1-simplex $\sig: \unin \to X$ that traverses each of the infinitely many circles in this space precisely once. Let $\beta\in H_1(X)$ denote the homology class of the 1-cycle $1\sig$, and note that \fe\ 1-cycle $\chi\in \beta$ there holds $\ell(\chi)=\infty$ because of the perpendicular segments. It is not hard to see that for $C':= \llbracket \beta \rrbracket \in \hompx$ we have $\inf \{\ell(\beta) \mid \beta \in C'\}=\infty$. Now let $\tau_i$ be a \clex\ that travels once around the circle of length $2^{-i}$ in $X$, and let $\psi_i$ denote the 1-chain $\sum_{j\leq i} \tau_j$. By \Lr{U} we can, \fe\ $i$, `patch' all circles of $X$ of length less than $2^{-i}$ to obtain an \arex\ $X_i$ of $X$ of some \evol\ $v(i)< \infty$ in which the 1-cycles $1\sig$ and $\psi_i$ are homologous. Note that $\lim_i v(i)=0$, thus $(\ec{\psi_i})_{\iin}$ is a \Cs\ equivalent to the constant sequence $(C')_{\iin}$, which means that $(1\tau_i)_{\iin}$ is a \rep\ of the class $C\in \lhomx$ containing these sequences. Thus $\ell(C)\leq  \sum_i \ell(\tau_i) = 1$.
\noproof
}

\showFig{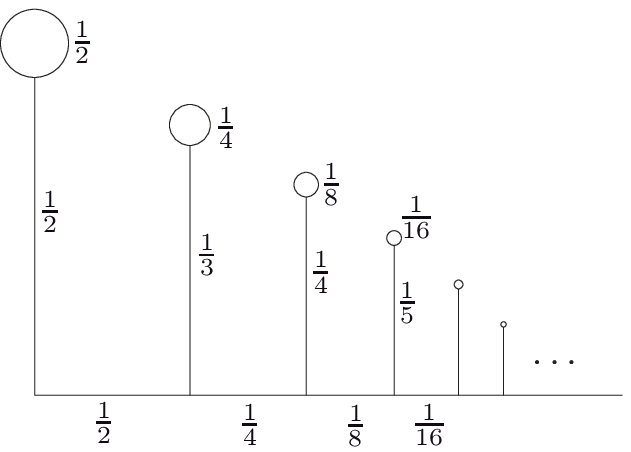}{A compact subspace of the real plane. The numbers denote the lengths of the corresponding segments.}

Finally, it is worth mentioning that we cannot relax the assertion of \Tr{lCycD} to require that $X$ is just complete rather than compact. For example, the cylinder $\{(x,y,z) \in \R^3 \mid z\geq 1, x^2 + y^2 = 1+ 1/z\}$ with the Euclidean metric is complete, but it is easy to see that no non-trivial element of \lhom\ has a \rep\ of minimum length. 

\section{Sketch of the main proof} \label{Ssketch}

The proof of our main result, \Tr{lCycD}, consists of two major steps: the first step is algebraic, and shows that every $C\in \lhom$ can be `decomposed' as a sum $\sum D_i$ of simpler elements of \lhom, called \defi{\prim} elements, that are easier to work with. The second step is more geometric, and proves the assertion for these primitive elements. 

Our intuition behind a \prim\ element is that it is a homology class corresponding to a single circle, and indeed we will prove, in \Sr{Sprim}, that every \prim\ element $D$ has a representative consisting of a \clex\ $z$, and in fact one of the desired length $\ell(z)=\ell(D)$. We obtain $z$ by a geometric construction: starting from a sequence of closed 1-simplices $\sig_i$ representing $D$ whose lengths converge to $\ell(D)$, we exploit the compactness of our space to find a subsequence that converges pointwise to the desired 1-simplex $z$, and show that $\ec{z}=D$  by constructing arbitrarily small metric cylinders joining $z$ to some $\sig_i$. See also \Er{excamera}, where we could choose $\tau$ to be the desired \clex\ $z$.

Now having a decomposition $C=\sum D_i$ as above, we can try to combine all the \clex es $z_i$ we got as representatives of each $D_i$ to form a \rep\ of $C$. But will such a \rep\ have the desired total length $\sum \ell(z_i) = \ell(C)$? In general not, if our decomposition is arbitrary. For example, in the graph of \fig{owl} consider the class $C = \llbracket \sig + \tau \rrbracket$. We could write $C = D_1 + D_2$ where $D_1 = \llbracket\sig \rrbracket$ and $D_2 = \llbracket\tau \rrbracket$ are both primitive. Now $\sig,\tau$ are \clex es that do attain the length of $D_1,D_2$ respectively, but we cannot combine them into a representative of $C$ of minimum length, because $\ell(\sig)+\ell(\tau)> \ell(C)$; indeed, $C$ has the representative $\rho$ whose length is smaller than $\ell(\sig)+\ell(\tau)$ because it avoids the middle edge. This example shows that if we want to follow the above plan of first decomposing $C$ as a sum of \prim\ elements and then combine shortest representatives of those elements into a \rep\ of $C$ of the desired total length $\ell(C)$, then our decomposition has to be `economical'. If our space is a graph then it is easy to say what `economical' should mean: no edge should be used in more than one summands. In a general space this is less obvious, but there is an elegant way around it described in \Sr{SSsplit}. We will prove, in \Sr{Ssplpr}, that every $C\in \lhom$ can be decomposed as a sum $\sum D_i$ where, not only the $D_i$ are \prim, but also the decomposition is economical in this sense. This proof is  algebraic, and we obtain a more general abstraction described in the next section. 

\showFig{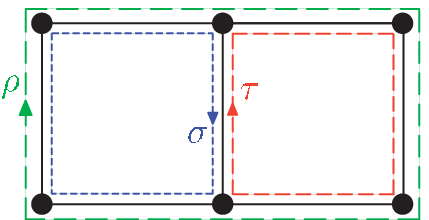}{A simple example showing that we need our primitive decompositions to be economical.}


\section{Intermezzo: generalising to abelian metrizable topological groups} \label{Smezzo}
In this section we state an intermediate result, mentioned also in \Sr{IntMain}, that might be useful in other contexts too. It says that if a topological group $H$ and an assignment $\ell: H \to \R^+$ (which can be thought of as an assignment of lengths) satisfy certain axioms, then every element of $H$ can be written as a sum of \defi{primitive} elements, which we define below, and this sum is in a sense `economical' (recall the discussion in the previous section). 

The reader will lose nothing by assuming that $H= \lhomx$ throughout this section. 

Given two elements $C,D$ of $H$ we will write $C\hleq D$ if $ \ell(C) = \ell(D) + \ell(C - D)$. Note that
\labtequ{compleq}{if $D\hleq C$ then $C-D\hleq C$,}
since $\ell(C - (C-D))= \ell(D)$.

We will say that an element $C$ of $H$ is \defi{\prim} if $C\neq \zero$ and \ti\ no $D\in H \sm \{C,\zero\}$ such that $D \hleq C$ holds. 

The reader may choose to skip the rest of this section, since this is a corollary of our main result rather than something that we will need later.

\begin{theorem} \label{Tmezzo}
Let $(H,+)$ be an abelian metrizable topological group, and suppose a function $\ell: H \to \R^+$ is given satisfying the following properties:

\begin{enumerate}
\item $\ell(C)=0$ \iff\ $C=\zero$; 
\item $\ell(C+D) \leq \ell(C) + \ell(D)$ \fe\ $C,D\in H$;
\item if $D = \lim C_i$ then $\ell(D) \leq \liminf \ell(C_i)$;
\item for some metric $d$ of $H$ \ti\ a bound $U\in \R$ \st\ $d(C,\zero) \leq U\ell^2(C)$ \fe\ $C\in H$ (i.e.\ an isoperimetric inequality holds). \label{piv}
\end{enumerate}

Then every element $C$ of $H$ can be represented as a (possibly infinite) sum $C=\sum D_i$ of \defi{primitive} elements $D_i$ so that $\ell(C)= \sum \ell(D_i)$. 
\end{theorem}

Infinite sums as in the conclusion of the theorem are formalised, in \Sr{Infsums}, using the concept of nets.

Since it is the group \lhom\ we are interested in in this paper, we will give a formal proof of \Tr{Tmezzo} only for the special case when $H=\lhom$ (more precisely, when $H$ is the subgroup of elements of \lhom\ with finite length). In this case \Tr{Tmezzo} is tantamount to  \Cr{split} below. However, the reader interested in \Tr{Tmezzo} in its full generality will easily be able to check that the same proof applies, since no other properties of \lhom\ are used in the proof of \Cr{split} than the conditions of \Tr{Tmezzo}.

One can relax condition \ref{piv} above a bit by replacing it with the following
\begin{equation*}  	\begin{minipage}[c]{0.9\textwidth}  
$(iv')$ 	if $C\in H$ is fragmentable then $C=\zero$.
 \end{minipage} \ignorespacesafterend \end{equation*}
The term \defi{fragmentable} is defined in \Sr{Eisop} below.

\section{Splitting homology classes into primitive subclasses} \label{Ssplpr}

The main result of this section, \Cr{split}, is that every $C\in \lhom$ can be written as a sum of primitive elements $D_i \hleq C$. This is the first step of the proof of our main result as sketched in \Sr{Ssketch}. Recall the definitions of \defi{primitive} and $\hleq$ from \Sr{Smezzo}.

\subsection{Infinite sums in \lhom} \label{Infsums}

For the proof of our main result we are going to use some standard machinery related to \defi{nets} in order to be able to rigorously define sums of infinitely many elements of \lhom. Let us first recall the necessary definitions. 

A \defi{net} in a topological space $Y$ is a function from some directed set $A$ to $Y$. A \defi{directed set} is a nonempty set $A$ together with a reflexive and transitive binary relation, that is, a preorder, with the additional property that every pair of elements has an upper bound in $A$. One can think of a net as a generalisation of the concept of a sequence, and one is usually interested in the convergence of such a generalised sequence: we say that the net $(x_\alpha)$ \defi{converges} to the point $y\in Y$, if for every neighborhood $U$ of $y$ \ti\ a $\beta\in A$ \st\ $x_\alpha\in U$ \fe\ $\alp\geq \beta$. See \cite{willard} for more details. In our case, the topological space $Y$ in which our nets will take their values will always be \lhom, bearing the topology induced by the metric $d_1$.

We will say that an infinite family \fml{C}\ of elements of \lhom\ is \defi{\unsu} if \fe\ $\eps>0$ \ti\ a finite subset $F$ of $\ci$ \sot\ \fe\ two finite sets $A,B\supseteq F$ \tho\ $d(\sum_{i\in A} C_i, \sum_{i\in B} C_i) < \eps$; in other words, if the family $\{\sum_{i\in F} C_i\}_{F\in \cf}$ is a Cauchy net, where $\cf$ is the set of finite subsets of $\ci$ preordered by the inclusion relation. Since \lhom\ is complete, it is well known that if  \fml{C}\ is \unsu\ then the net $\{\sum_{i\in F} C_i\}_{F\in \cf}$ converges to an element $C\in \lhom$, see \cite[Proposition 2.1.49]{megginson}. In this case, we call $C$ the \defi{sum} of the family \fml{C}\ and write $C = \sum_{i\in\ci} C_i$. Note that if $\ci$ is countable then \fe\ enumeration $a_1, a_2, \ldots$ of $\ci$ \tho\ 
\labtequ{ctblsum}{$\sum_{i\in\ci} C_i = \lim_i \sum_{1\leq j \leq i} C_{a_j}$.}

Our next lemma generalises the triangle inequality for \lhom\ (\Cr{ltriang}) to infinite sums using the notions we just defined.

\begin{lemma}\label{infsum}	
Let \fml{C}\ be an \unsu\ family of elements of $\lhom$. Then $\ell(\sum_{i\in\ci} C_i) \leq \sum_{i\in\ci} \ell(C_i) $.
\end{lemma}
\begin{proof}

If $\sum_{i\in\ci} \ell(C_i)=\infty $ then there is nothing to show, so suppose that $\sum_{i\in\ci} \ell(C_i) < \infty $.

We may assume \obda\ that $\ell(C_i)>0$ holds \fe\ $i\in\ci$, for if $\ell(C_i)=0$ then $C_i= \zero$ by \Or{zerol}. Thus, we may also assume that $\ci$ is countable, and let $a_1,a_2,\ldots$ be an enumeration of $\ci$. 

Let $C:=\sum_{i\in\ci} C_i$ and let $D_i:= \sum_{1\leq j \leq i} C_{a_j}$ \fe\ $i$. We have $ C = \lim_i D_i$ by \eqref{ctblsum}. By the definition of $\ell(C)$ we then have 
\labtequ{CD}{$\ell(C)\leq \lim \ell(D_i)$.}
Applying \Cr{ltriang} (several times) to $D_i$ we obtain $\ell(D_i) \leq \sum_{1\leq j \leq i} \ell(C_{a_j}) < \sum_{i\in\ci} \ell(C_i)$. Combining this with \eqref{CD} yields $\ell(C)\leq \sum_{i\in\ci} \ell(C_i)$ as desired.



\end{proof}

\subsection{Splitting homology classes into shorter ones} \label{SSsplit}

We introduce the notation $C = D \blcup E$ to denote the assertion that $C = D + E$ and $\ell(C) = \ell(D) + \ell(E)$. Note that this definition implements the intuition outlined in \Sr{Ssketch} that $D + E$ is an `economical' way to split $C$.
It follows by the definitions  that 
\labtequ{compleqI}{$D\hleq C$ \iff\ $C = D \blcup (C-D)$.} 

More generally, the notation $C = \blcup_{i\in \ck} D_i $ (or $D_1 \blcup \ldots \blcup D_k$), where $\ck$ is a possibly infinite set of indices, denotes the assertion that $C = \sum_{i\in \ck} D_i$ and $\ell(C) = \sum_{i\in \ck} \ell(D_i)$.

Our next lemma shows that, in a sense, $\blcup$ behaves well \wrt\ composition:

\begin{lemma}\label{trans}
Let $C,D,E,F,G \in \lhom$ be such that $C =D \blcup E$ and $E =F \blcup G$. Then the following assertions hold:
\begin{enumerate}
\item \label{ti} $C =D \blcup F \blcup G$; 
\item \label{tii} $C = (D + F) \blcup G$, and
\item \label{tiii} $D + F = D \blcup F$. 
\end{enumerate}
In particular, $F, G, (D+F) \hleq C$.
\end{lemma} 
\begin{proof}
By the assumptions we have $\ell(C) = \ell(D) + \ell(E) =  \ell(D) + \ell(F) + \ell(G)$, which yields \ref{ti}.

Note that $C = D + F + G$. By \Cr{ltriang} we have $\ell(C) = \ell( (D + F) + G ) \leq \ell(D + F) + \ell( G )$, and $\ell(D + F)\leq \ell(D) + \ell(F)$. Now since we have already proved that $\ell(C) = \ell(D) + \ell(F) + \ell( G )$, equality must hold in both above inequalities. The first of these equalities yields \ref{tii} and the second yields \ref{tiii}.
\end{proof}

\comment{
	\begin{lemma}\label{trans}
	For every $C,D,E \in \lhom$, if $E \hleq D \hleq C$ then $E \hleq C$.
\end{lemma} 
\begin{proof}
By \Cr{ltriang} we have $$\ell(C - E) = \ell( (C-D) + (D - E) ) \leq \ell(C-D) +  \ell(D-E).$$ By the assumption we have $\ell(C - D) = \ell(C) - \ell(D)$ and $\ell(D-E) = \ell(D) - \ell(E)$. Replacing in the above inequality we obtain $$\ell(C-E) \leq \ell(C) - \ell(D) + \ell(D) - \ell(E) = \ell(C)- \ell(E).$$ On the other hand, applying \Cr{ltriang} again we obtain $\ell(C-E) \geq \ell(C)- \ell(E)$. The latter two inequalities together yield $\ell(C-E) = \ell(C)- \ell(E)$ i.e.\ $E\hleq C$.
	\end{proof}
}

This nice behaviour of $\blcup$ extends to infinite sums too:

\begin{lemma}\label{lemPlus}
If $\ell(C)<\infty$ and $C = \blcup_{i\in \ck} D_i$ then \fe\ subset $\cm\subseteq \ck$ there holds $\sum_{i\in \cm} D_i =\blcup_{i\in \cm} D_i$ and $C = \sum_{i\in \cm} D_i \blcup \sum_{i\in \overline{\cm}} D_i$, where $\overline{\cm}:= \ck \sm \cm$. 
\end{lemma} 
\begin{proof}
By \Cr{ltriang} we have $\ell(C) \leq \ell(\sum_{i\in \cm} D_i) + \ell(\sum_{i\in \overline{\cm}} D_i)$ and by \Lr{infsum} we have $\ell(\sum_{i\in \overline{\cm}} D_i) \leq \sum_{i\in \overline{\cm}} \ell(D_i)$. Combining the last two inequalities we obtain
\begin{eqnarray*}
\ell(\sum_{i\in \cm} D_i) \geq \ell(C) - \ell(\sum_{i\in \overline{\cm}} D_i) = &\sum_{i\in \ck} \ell(D_i) &- \ell(\sum_{i\in \overline{\cm}} D_i) \geq\\
&\sum_{i\in \ck} \ell(D_i) &- \sum_{i\in \overline{\cm}} \ell(D_i) =
\sum_{i\in \cm} \ell(D_i),
\end{eqnarray*}
where we used our assumption that $\sum_{i\in \ck} \ell(D_i) =\ell(C) <\infty$. 
Applying \Lr{infsum} again we also have 
$$\ell(\sum_{i\in \cm} D_i) \leq \sum_{i\in \cm} \ell(D_i),$$
hence equality holds in the last two inequalities, which proves that 
$$\sum_{i\in \cm} D_i =\blcup_{i\in \cm} D_i.$$ 

Similarly, we have  $\sum_{i\in \overline{\cm}} D_i =\blcup_{i\in \overline{\cm}} D_i$. The assertion 
$$C = \sum_{i\in \cm} D_i \blcup \sum_{i\in \overline{\cm}} D_i$$
now easily follows from the definitions.
\end{proof}

\subsection{Exploiting the isoperimetric inequality} \label{Eisop}

We will say that an element $C\in \lhom$ is \defi{\del-fragmentable}, for some $\del\in \R_+$, if \ti\ a finite family \ffml{D}{\ck}, $D_i\in \lhom$, such that $C= \blcup_{i\in \ck} D_i$
and \fe\ $i$ there holds $\ell(D_i)< \del$. We will call $C$ \defi{fragmentable} if it is \del-fragmentable for arbitrarily small \del. It turns out that the only fragmentable element of \lhom\ is $\zero$:

\begin{lemma}\label{fragm}	
If $C\in\lhom$ is fragmentable then $C=\zero$.
\end{lemma}
\begin{proof}
Suppose $C$ is fragmentable, and fix some $\eps>0$ for which we want to show that $d_1( C, \zero)< \eps$. Let \ffml{D}{\ck}\ be a family witnessing the fact that $C$ is \del-fragmentable for some parameter \del\ that  we will specify later.

For every $\iin$,  we can,  by the definition of $\ell(D_i)$, find elements of \shom\ arbitrarily close (with respect to $d_1$) to $D_i$ whose lengths are arbitrarily close to $\ell(D_i)$; more formally, it follows by the definitions that \ti\ a class $\alp_i\in \shom$ with $\ell(\alp_i)< \ell(D_i) + \min(\del, \eps/2|\ck|)$ \st\ 
\labtequ{dbet}{$d_1( [\alp_i] , D_i ) < \eps/2|\ck|$.}


By the definition of $\ell(\alp_i)$ there is an 1-chain $\chi_i\in \alp_i$ \st\ $\ell(\chi_i) < \ell(\alp_i) + \min(\del, \eps/2|\ck|)$. Combining this with our assumption that $\ell(D_i)<\del $ and the choice of $\alp_i$ 
we obtain 
\labtequ{triadel}{$\ell(\chi_i) < \ell(\alp_i) + \del <  \ell(D_i) + 2\del < 3\del.$} 

By \Lr{U} \ti\ an extension $X_i$ of $X$ of \evol\ at most $U\ell^2(\chi_i)$ in which $\chi_i$ is null-homologous. Combining these extensions $X_i$ \fe\ $i$ we obtain an extension $X'$ of $X$ of \evol\ $V$ at most $U \sum_{i\in \ck} \ell^2(\chi_i)$ in which the 1-chain $\sum_{i\in \ck} \chi_i$ is null-homologous. Note that by the choice of the $\chi_i$ and the $\alp_i$ we have
$$\sum \ell(\chi_i) < \sum (\ell(\alp_i) + \eps /2|\ck|) < \sum (\ell(D_i) + \eps /2|\ck|+\eps /2|\ck|) = (\sum \ell(D_i))  + \eps = \ell(C)  + \eps,$$ 
where we used our assumption that $\ell(C)= \sum \ell(D_i)$. This means that $\sum \ell(\chi_i)$ is bounded from above; thus by~\eqref{squares} and \eqref{triadel} choosing \del\ small enough we can make $V$ arbitrarily small; in particular, we could have chosen a \del\ for which $V< \eps/2$ holds, which would imply 
\labtequ{chiz}{$d_1(\ec{\sum \chi_i} , \zero) < \eps/2$.}
Since $C = \sum D_i$ we easily obtain by \Cr{ltriang} and~\eqref{dbet} 
$$d_1( C,  \ec{\sum_{i\in \ck} \chi_i}) = d_1( \sum_{i\in \ck} D_i,  \ec{\sum_{i\in \ck} \chi_i}) \leq \sum_{i\in \ck} d_1( D_i, \ec{\chi_i} ) \leq \sum_{i\in \ck} \eps/2|\ck| = \eps/2,$$
and combined with \eqref{chiz} this yields $d_1( C, \zero) < \eps$, and proves that $C = \zero$ in this case.
\end{proof}

\subsection{A technical lemma}

The following somewhat technical lemma will be used in the proof of the main result of this section; it allows us to prove, using Zorn's Lemma, the existence of maximal families with certain properties. 

\begin{lemma}\label{infami}
Let \ofml{D}{\gam}\ be a family of elements $D_\alp$ of $\lhom\sm \{\zero\}$, indexed by an ordinal \gam,
\st\  \fe\ $\bet<\gam$ there holds $\sum_{\alp\leq \bet} D_\alp = \blcup_{\alp\leq \bet} D_\alp$ and $\sum_{\alp\leq \bet} D_\alp \hleq C$ for some fixed $C\in \lhom$ with $\ell(C)< \infty$.
Then \oseq{D}{\gam}\ is \unsu\ and \tho\ $\sum_{\alp < \gam} D_\alp \hleq C$ and $\sum_{\alp < \gam} D_\alp =\blcup_{\alp < \gam} D_\alp$.
\end{lemma} 
\begin{proof}

Since $D_\alp \neq \zero$, \Or{zerol} implies that $\ell(D_\alp)>0$ \fe\ $\alp<\gam$. 
As we are assuming that 
$\ell(\sum_{\alp\leq \bet} D_\alp) = \sum_{\alp\leq \bet} \ell(D_\alp)$ and that $\ell(\sum_{\alp\leq \bet} D_\alp)\leq \ell(C)$ \fe\ $\bet \le\gam$,  we have $\sum_{\alp\leq \bet} \ell(D_\alp) \leq \ell(C)<\infty$ \fe\ $\bet <\gam$, which implies that $\gam$ is countable and 
\labtequ{gaminfty}{$\sum_{\alp<\gam} \ell(D_\alp) \leq \ell(C) < \infty$.} 

Let $a_1,a_2,\ldots$ be an enumeration of $\gam$.
To see that \oseq{D}{\gam}\ is \unsu, note that \fe\ $\eps>0$ \ti\ an $\nin$ \st\ \fe\ $k>n$ \tho\ $\ell(\sum_{n< j < k} D_{a_j} ) \leq \sum_{n< j < k} \ell( D_{a_j} ) < \eps$, hence $d_1(\sum_{n< j < k} D_{a_j}, \zero) <U\eps^2$  by \Lr{U} and the definition of $d_1$.

Thus $S:= \sum_{\alp<\gam} D_\alp$ is well-defined (see \Sr{Infsums}), and by \eqref{ctblsum} we have 
\labtequ{limlim}{$S = \lim_n D^n$ and $C - S = \lim_n (C - D^n)$,}
where $D^n:= \sum_{1\leq j\leq n} D_{a_j}$. 

We have to prove that $S \hleq C$, i.e.\ that $\ell(C) = \ell(S) + \ell(C-S)$. By the definition of $\ell()$ and \eqref{limlim} we obtain $\ell(S) \leq \liminf_n \ell(D^n)$ and $\ell(C - S) \leq \liminf_n \ell(C -D^n)$. However, \eqref{gaminfty} and \Cr{ltriang} easily imply that both $(\ell(D^n))$ and $(\ell(C-D^n))$ converge, and so we can write
\labtequ{SCS}{$\ell(S) \leq \lim_n \ell(D^n)$ and $\ell(C - S) \leq \lim_n \ell(C -D^n)$.}
Combining this with the fact that $\ell(C) \leq \ell(S) + \ell(C-S)$, which we obtain from  \Cr{ltriang}, we have
\labtequ{CDi}{$\ell(C) \leq \lim_n \ell(D^n) + \lim_n \ell(C -D^n) = \lim_n (\ell(D^n) + \ell(C -D^n))$.}

We claim that $D^n\hleq C$ holds \fe\ $n$. Indeed, note that \ti\ a $\bet<\gam$ \st\ $a_j\leq \bet$ holds \fe\ $j\leq n$. As we are assuming that $\sum_{\alp\leq \bet} D_\alp = \blcup_{\alp\leq \bet} D_\alp$, \Lr{lemPlus} (for the application of which we set $C = \sum_{\alp\leq \bet} D_\alp$, $\ck=\bet$ and $\cm=\{a_1,\ldots, a_n\}$) implies that 
\labtequ{Dn}{$D^n = \blcup_{1\leq j\leq n} D_{a_j}$}
and that $D^n\hleq \sum_{\alp\leq \bet} D_\alp$.
As we are furthermore assuming that $\sum_{\alp\leq \bet} D_\alp \hleq C$ holds, the transitivity of $\hleq$ (see \Lr{trans}) implies that $D^n\hleq C$ as claimed. 

This means that $\ell(D^n) + \ell(C -D^n) = \ell(C)$ \fe\ $\nin$. Plugging this into \eqref{CDi} yields $\ell(C) \leq \lim_n \ell(C) = \ell(C)$. Thus equality must hold throughout in \eqref{SCS} and \eqref{CDi}. This implies that $\ell(C) = \ell(S) + \ell(C-S)$ ---i.e.\  $\sum_{\alp<\gam} D_\alp \hleq C$--- and that $\ell(S) = \lim_n \ell(D^n)$. Using \eqref{Dn} the latter yields $\ell(S) = \ell(\sum_{\alp<\gam} D_\alp) =\lim_n \sum_{1\leq j\leq n} \ell(D_{a_j}) =\sum_{\alp<\gam} \ell(D_\alp)$ as desired.

\comment{
		We claim that 
	\labtequ{ccauchy}{\Fe\ $\eps'>0$ \ti\ a $\bet< \gam$ \st\ $d_1(D_\alp, D_\alp')< \eps'$ holds \fe\ $\alp,\alp'>\bet$.}

$\ell(C_\bet) \leq \ell(C)$ \fe\ $\bet<\gam$, 

$\ell(C_i) \leq \ell(C_j)$ \fe\ $i< j \in\N$

thus \ti\ $m\in\N$ \st\ $\ell(C_i) - \ell(C_j) < \eps$ \fe\ $i,j>m$

Since $C_i\hleq C_j$ we have $\ell(C_j - C_i)= \ell(C_j) - \ell(C_i) < \eps$

By \Cr{obl}\ \ti\ an 1-cycle $\phi$ \st\ $\ell(\phi)< \ell(C_j - C_i) + \eps$ and $C_1( (C_j - C_i),  \ec{\phi} ) < \eps$.

Thus by \Lr{U} 
$$d_1( (C_j - C_i) , \zero) \leq d_1( (C_j - C_i) , \ec{\phi} ) + d_1( \ec{\phi} , \zero) \leq  \eps + U(2\eps)^2.$$

		It follows easily from the definition of $d_n$ that $d_1( C_j, C_i) \leq d_1( (C_j - C_i) , \zero)$, and with the previous inequality we obtain $d_1( C_j, C_i) \leq \eps'$ where $\eps':= \eps + U(2\eps)^2$. As $U$ is fixed we can make $\eps'$ arbitrarily small by choosing $\eps$ small enough, hence 		claim~\eqref{ccauchy} follows.
}
\end{proof}

\subsection{Existence of primitive decompositions}

We can now complete the proof of the main result of this section, that every non-trivial element $C$ of $\lhomx$ can be decomposed as a sum of primitive elements. We do this by first proving that we can find at least one primitive element in $C$, and then using Zorn's Lemma to find a maximal family of primitive elements in $C$. The former task is fulfilled by the following lemma.

\begin{lemma}\label{prim}
For every $C\neq \zero\in \lhomx$ with $\ell(C)<\infty$ there is $D\hleq C$ \st\ $D$ is \prim.
\end{lemma} 
\begin{proof}
By \Lr{fragm}, $C$ is not \lam-fragmentable for some $\lam\in \R_+$, i.e.\ for every finite family \ffml{D}{\ck}\ such that $C= \blcup_{i\in \ck} D_i$ there is a member $D_j$ with $\ell(D_j)\geq \lam$. Note that there is also no infinite family \ffml{D'}{\ci}\ such that $C= \blcup_{i\in \ci} D'_i$ and \fe\ $i$ there holds $\ell(D_i)< \lam$. For if such a family exists, then we can find a finite subfamily \ffml{D'}{\ci'}\ such that $\sum_{i\in \ci'} \ell(D'_i) > \ell(C) - \lam$, which implies $\sum_{i\in \ci\sm\ci'} \ell(D'_i) < \lam$, and hence $\ell(\sum_{i\in \ci\sm\ci'} D'_i) < \lam$ by \Lr{infsum}. But then, extending \ffml{D'}{\ci'}\ by one member, namely $\sum_{i\in \ci\sm\ci'} D'_i$, we obtain a finite family \ffml{D'}{\ci''} which satisfies $\ell(D'_i) <\lam$ \fe\ $i\in \ci''$, and it is not hard to see that $C= \blcup_{i\in \ci''} D'_i$ holds; this contradicts the fact that $C$ is not \lam-fragmentable.


\comment{
	Now let \ffml{D}{\gamma}, $D_i\in \lhom\sm \zero$ be a family, indexed by an ordinal number $\gamma$, that is maximal with the following properties:
\begin{enumerate}
\item \label{ri} $D_\beta \hleq C - \sum_{\alpha < \beta} D_\alp$ \fe\ $\beta<\gam$; in particular, $D_0 \hleq C$;
\item \label{rii} $\sum_{\alp\in\gam} \ell(D_\alp) \leq \ell(C) - \lam$, and
	\item \label{riii} $D:= \sum_{\alp\in\gam} D_\alp \hleq C$.
	\end{enumerate}
}

Now let \ofml{D}{\gamma}, $D_\alp\in \lhom\sm \zero$ be an \unsu\ family, indexed by an ordinal number $\gamma$, that is maximal with the following properties:
\begin{enumerate}
\item \label{ri} $\sum_{\alpha < \beta} D_\alp \hleq C$ \fe\ $\beta\leq\gam$;
\item \label{rih} $\sum_{\alpha < \beta} D_\alp = \blcup_{\alpha < \beta} D_\alp$ \fe\ $\beta\leq\gam$, and 
\item \label{rii} $\sum_{\alp<\gam} \ell(D_\alp) \leq \ell(C) - \lam$.
\end{enumerate}
To see that a maximal such family exists, apply Zorn's Lemma on the set of all such families ordered by the subfamily relation, using \Lr{infami} in order to show that every chain has an upper bound. We are not yet assuming that this maximal family is non-trivial.

Let $D:= \sum_{\alp<\gam} D_\alp$, and note that $D\hleq C$ by \ref{ri}. It is not hard to see that either $D$ or $C-D$ (or both) is still not \lam-fragmentable, for if both split into families with elements of lengths less that \lam, then so does $C$; more formally, suppose there are finite families \ffml{D}{\cm}\ and \ffml{D}{\cn}\ such that $D= \blcup_{i\in \cm} D_i$, $C-D= \blcup_{i\in \cn} D_i$, and $\ell(D_i)<\lam$ \fe\ $i\in \cm \cup \cn$. We claim that $C = \blcup_{i\in \cm \cup \cn} D_i$. Easily, $C = \sum_{i\in \cm \cup \cn} D_i$. To see that $\ell(C) = \sum_{i\in \cm \cup \cn} \ell(D_i)$, recall that $\ell(C) = \ell(D) + \ell(C-D)$ by \ref{ri}, that  $\ell(D) = \sum_{i\in \cm} \ell(D_i)$, and that $\ell(C-D) = \sum_{i\in \cn} \ell(D_i)$. This proves that either $D$ or $C-D$ is not \lam-fragmentable.

But if $C-D$ is not \lam-fragmentable, then it is primitive: for if \ti\ an $F_0\neq \zero$ with $F_0\hleq C-D$ and $F_0\neq C-D$, then either $F_0$ or $F_1:= C - D - F_0$ has length at least \lam\ since $C -D$ is not \lam-fragmentable and,  by \eqref{compleqI}, $C - D = F_0 \blcup F_1$. Assume \obda\ that $\ell(F_1) \geq \lam$; we can now enlarge the family \ofml{D}{\gam}\ by one member, namely $F_0$, to obtain a new family \ofml{D}{\gam^+}\ that  contradicts the maximality of \fml{D}: to prove that \ofml{D}{\gam^+}\ also satisfies requirement~\ref{ri} it suffices to check that $\sum_{\alpha < \gam^+} D_\alp \hleq C$. We have $\sum_{\alpha < \gam^+} D_\alp = D + F_0$ by construction, and by assertion~\ref{tii} of \Lr{trans} we obtain
\labtequ{DF}{$D + F_0 \hleq C$,}
which proves that \ofml{D}{\gam^+}\ satisfies~\ref{ri}. To prove that \ofml{D}{\gam^+}\ also satisfies requirement~\ref{rih}, it suffices again to consider the case $\bet=\gam^+$; in other words, to prove that $D + F_0 = \blcup_{\alpha < \gam^+} D_\alp$. Thus we have to prove that $\ell(D + F_0 ) = \sum_{\alpha < \gam^+} \ell(D_\alp) = \ell(D) + \ell(F_0)$, where for the last equality we used the fact that \ref{rih} holds for $\bet=\gam$ and $D_{\gam}= F_0$. But this follows from assertion~\ref{tiii} of \Lr{trans},
and so \ofml{D}{\gam^+}\ also satisfies~\ref{rih}. Finally, to see that \ofml{D}{\gam^+}\ satisfies~\ref{rii}, note that $\sum_{\alp\in\gam^+} D_\alp = D + F_0$, that $\ell(C) = \ell(D+ F_0) + \ell( F_1 )$ by \eqref{DF}, and that $\ell( F_1 ) \geq \lam$. This completes the proof that if $C-D$ is not \lam-fragmentable then $C-D$ is primitive, for otherwise the maximality of \ofml{D}{\gam} is contradicted.


Thus, if $C-D$ is not \lam-fragmentable then we are done, since $D\hleq C$ and so we also have $C-D\hleq C$ by \eqref{compleq}. So suppose it is not, in which case it is $D$ that is  not \lam-fragmentable. Recall that $\ell(D)\leq \ell(C) - \lam$ by~\ref{rii}.

To sum up, having assumed that $C$ is not \lam-fragmentable, we proved that either there is a primitive $B \hleq C$, in which case we are done, or there is a $D \hleq C$ that is also not \lam-fragmentable (for the same \lam) and satisfies $\ell(D)\leq \ell(C) - \lam$. In the latter case, we can repeat the whole argument replacing $C$ with $C_1:= D$; this will again yield either a primitive $B \hleq C_1$, or a $C_2 \hleq C_1 $ that is also not \lam-fragmentable and satisfies $\ell(C_2)\leq \ell(C_1) - \lam \leq \ell(C) - 2\lam$, and so on. But as $\ell(C)$ is finite and \lam\ positive, this procedure must stop after finitely many steps, yielding a primitive $B \hleq C_j \hleq C_{j-1} \ldots \hleq C$. As $\hleq$ is transitive (\Lr{trans}) we obtain $B \hleq C$. This completes the proof.
\end{proof}

We can now state and prove the main result of this section.

\begin{corollary}\label{split}
For every $C\neq 0\in \lhomx$ with $\ell(C)<\infty$ there is a family \fml{D}\ of \prim\ elements of \lhomx\ \st\ $C= \blcup_{i\in I} D_i$.
\end{corollary}
\begin{proof} 	
Using Zorn's Lemma we find a maximal family \ofml{D}{\gam}\ of primitive $D_\alp\in \lhom$ \st\ 
\begin{enumerate}
\item \label{si} $\sum_{\alpha < \beta} D_\alp \hleq C$ \fe\ $\beta\leq\gam$; and
\item \label{sii} $\sum_{\alpha < \beta} D_\alp = \blcup_{\alpha < \beta} D_\alp$ \fe\ $\beta\leq\gam$. 
\end{enumerate}
Indeed, consider the set of all such families ordered by the subfamily relation, and apply \Lr{infami} in order to show that every chain has an upper bound.

Let $D:=\blcup_{i\in I} D_i$. We claim that $C - D = 0$. For suppose not. Then  by \Lr{prim} there is a primitive $F\hleq C-D$. Now extend the family \ofml{D}{\gam}\ by one member $D_\gam:= F$. By \ref{tii} of \Lr{trans} the new family still satisfies \ref{si}. To prove that it also satisfies \ref{sii} we only have to show that $\ell(\sum_{\alpha \leq \gam} D_\alp) = \sum_{\alpha \leq \gam} \ell(D_\alp) = \ell(D) + \ell(F)$ (where we used the fact that the original family satisfies \ref{sii}), but this follows from \ref{tiii} of \Lr{trans}. Thus the extended family contradicts the maximality of \ofml{D}{\gam}, which proves that $C - D = 0$ and establishes our assertion.
\end{proof}

\section{Proof for primitive elements} \label{Sprim}

By \Cr{split} every non-trivial element $C$ of \lhom\ can be written as a sum of primitive elements $D_i$ so that $\ell(C) = \sum \ell(D_i)$. All that remains to show is that our main theorem holds for those elements $D_i$:

\begin{lemma}\label{lem}
If $D\in \lhom$ is primitive then there is a \clex\ $z$ \st\ $D = \ec{z}$ and $\ell(z)=\ell(D)$.
\end{lemma} 
\begin{proof} 
We are going to obtain the desired closed simplex $z$ as a limit, in a sense, of a sequence of closed simplices $\sig^1_i$ related to $D$. Our proof is organised in three steps. In the first step we construct this sequence $(\sig^1_i)$. In the second step we construct $z$ and, at the same time,  homotopies between $z$ and the $\sig^1_i$ in appropriate \arex s of $X$, implying that $\ec{z} = \lim \ec{sig^1_i}$. Finally, in a third step we show that $D = \ec{z}$ and that $\ell(z)= \ell(D)$. We then remark that the closed simplex $z$ we constructed must indeed be a \clex.

\subsection*{Step I: the sequence $(\sig^1_i)$}

By \Cr{obl} \ti\ a sequence of  1-cycles \seq{\chi} \st\ $(\ec{\chi_i})_{\iin}$ is a \Cs\ in $D$ and  $\ell(D) = \lim_i \ell(\chi_i)$.

By concatenating some of the simplices in $\chi_i$ if necessary, we may assume \obda\ that every simplex in $\chi_i$ is closed. For every $i$ enumerate the (closed) simplices in $\chi_i$ as $\sig_i^1, \ldots, \sig_i^{k_i}$ in such a way that 
\labtequ{ord}{$\ell(\sig_i^j) \geq \ell(\sig_i^m)$ if $j < m$.}
For convenience, if $m>k_i$ then we let $\sig_i^m$ denote a trivial 1-simplex in $X$ (thus $\ell(\sig_i^m)=0$ for $m>k_i$).

Let $\cm\subseteq \N$ be the set of superscripts $m$ \st\ \seq{\sig^m}\ has no infinite subsequence \susq{\sig^m}{\alpha}\ \st\ $\lim_i \ell(\sig^m_{\alpha_i}) = 0$. Note that, by \eqref{ord}, 
\labtequ{M}{if $m\in \cm$ then $\{{1,\ldots, m-1}\} \subset \cm$.}

We begin with a simple and instructive fact indicating the significance of \cm:
\begin{noclaim}
if $\cm = \emptyset$ then $D = \zero$.
\end{noclaim}

Indeed, if $\cm = \emptyset$ then there is an infinite subsequence \susq{\sig}{\alp}\ of \seq{\sig} \st\ $lim_i \ell(\sig^1_{\alpha_i}) = 0$. We will show that \fe\ $\eps>0$ there holds $d_1(D,\zero)<\eps$. For this, pick $j=\alp_k\in \N$ large enough that 
\begin{enumerate}
\item \label{ji} $d_1(D,\ec{\chi_j})<\eps/2$;
\item $\ell(\chi_j)< \ell(D) + \eps$, and
\item \label{jiii} $\ell(\sig^1_j) < \ell(D) \lam$,
\end{enumerate}
where $\lam=\lam(\eps)\in \R^+$ is some parameter that we will choose later. By \eqref{ord} we have $\ell(\sig^m_{j}) < \ell(D) \lam$ \fe\ $m\in \N$. By \Lr{U} \ti\ \fe\ $m$ an \arex\ $X_m$ of $X$ of
\evol\ at most $U\ell^2(\sig^m_j)$ in which $\sig^m_j$ is null-homologous. Combining all these \arex s we obtain a single \arex\ $X_\eps$ of $X$ of \evol\ at most $v:=\sum_{m \in \N} U\ell^2(\sig^m_j)$ in which $\chi_j$ is null-homologous. This means that 
\labtequ{chij}{$d_1(\ec{\chi_j},\zero)\leq v$.}

Since $\sum_{m \in \N} \ell(\sig^m_j) = \ell(\chi_j)< \ell(D) + \eps$, given $\ell(D)$ and \eps\ we can, by \eqref{squares} and \ref{jiii}, choose $\lam$ small enough that $v<\eps/2$. As $d_1(D,\zero) \leq d_1(D,\ec{\chi_j}) +d_1(\ec{\chi_j},\zero) < \eps/2 + v$ by \ref{ji} and~\eqref{chij}, and \eps\ was chosen arbitrarily, we have proved the Claim.

\medskip


As we are assuming that $D$ is primitive, the Claim proves that $\cm \neq \emptyset$, and thus $1\in M$ by \eqref{M}.

We may assume \obda\ that 
\labtequ{wp}{$\sig^1_i$ has \cosp\ \fe\ $i$.}

We are going to construct $z$ as a `limit' of the $\sig^1_i$ (it will turn out that $\cm=\{1\}$). For this, let \susq{\chi}{a}\ be a subsequence of \seq{\chi}\ \st\ $\lim_i \ell(\sig^1_{a_i})=: r$ exists. Note that we have already proved that $r>0$. Moreover, $r<\infty$ holds since $C$ is \prim\ and thus, easily, $\ell(C)<\infty$.

\note{
\begin{proposition}
Let $\seq{f}$ be a sequence of topological paths with \cosp\ in a fixed metric space such that $\ell(f_i)$ converges and $\lim \seq{f}=f$~pointwise. Then $f$ is continuous.
\end{proposition} 
\begin{proof}[Proof (sketch).]
Let $r:=\lim \ell(f_n)$. Consider $x\in \unin$ and an open ball $B_\eps(f(x))$ around $f(x)$ with $\eps<< r$. Let $n_0$ be large enough that $d( f(x), f_i(x) ) < \eps/2$ and $|r - \ell(f_i)| <\eps/2$ \fe\ $i>n_0$. Pick a (small) $\del\in\R$ and note that since the $f_i$ have \cosp,  $f_i\restr [x-\del,x+\del]$ is the concatenation of two topological paths of length at most $(r+\eps)\del$ \fe\ $i>n_0$. Thus $f_i([x-\del,x+\del])$ is contained in $B_{(r+\eps)\del}{f_i(x)}$ \fe\ $i>n_0$, which implies that $f([x-\del,x+\del])$ is contained in $B_{\eps}(f(x))$ if we choose \del\ small enough compared to $r$ and \eps. 
\end{proof}
}

It is not hard to see that there is a subsequence \susq{\sig^1}{b}\ of \susq{\sig^1}{a}\ \st\ the restrictions 
\labtequ{pwise}{$\sig^1_{b_i}\restr \Q$ converge pointwise.}
Indeed, let $q_1,q_2,\ldots$ be an enumeration of \Q. Find a subsequence \seq{\tau^0}\ of \susq{\sig^1}{a}\ \st\ the points $\tau^0_{i}(q_1)$ converge. Then find a subsequence \seq{\tau^1}\ of \seq{\tau^0}\ \st\ the points $\tau^1_{i}(q_2)$ also converge, and so on. Now letting $\sig^1_{b_i}= \tau^i_i$ satisfies \ref{pwise}.
(We could have chosen any dense countable subset of \unin\ instead of \Q.)

\subsection*{Step II: Construction of $z$ and $h$}

By \eqref{wp} and  \eqref{pwise} it follows easily that
\labtequ{del}{\fe\ \del\ \ti\ an $n\in\N$ \st\ $\sig^1_{b_i}$ and $\sig^1_{b_j}$ are \del-close \fe\ $i,j\geq n$.}

Using \Lr{Uc} and~\eqref{del} 
we can now construct a subsequence $\susq{\sig^1}{c}$ of $\susq{\sig^1}{b}$ \st\ \fe\ $i$ \ti\ an \arex\ $X'_i$ of $X$ of \evol\ at most $2^{-i}$ in which $\sig^1_{c_i}$ and $\sig^1_{c_{i+1}}$ are homotopic: for every $i=0,1,\ldots$, use~\eqref{del} to obtain a $c_i$ \st\  $\sig^1_{b_i}$ and $\sig^1_{b_j}$ are $f(r,2^{-i})$-close \fe\ $i,j\geq c_i$, where the function $f$ is that of \Lr{Uc}. Choosing $c_i$ larger if needed, we may also ensure that  $c_i> c_{i-1}$ (where we set $c_{-1}:=0$), and that $\ell(\sig^1_{b_i})< r + 2^{-i}$ \fe\ $i\geq c_i$. Then, by \Lr{Uc}, there is indeed an extension $X'_i$ as desired. Let $h_i$ be a homotopy from  $\sig^1_{c_{i+1}}$ to  $\sig^1_{c_i}$ in $X'_i$ as supplied by \Lr{Uc}.

Combining all $h_i$ together we can obtain a continuous function $h':  (0,1] \times \unin \to X'$, where $X':= \bigcup X'_i$. We are later going to ``complete'' $h'$ into a homotopy $h:\unin \times \unin \to X'$ \st\ $h(0,x)$ is our desired simplex $z$. To define $h'$, suppose that \fe\ $i$ we had chosen the domain of $h_i$ to be $[2^{-{i}},2^{-(i+1)}] \times \unin$. Intuitively, the interval $[2^{-{i}},2^{-{i+1}}]$ here corresponds to `time'; think of time as running in the negative direction if you prefer the homotopies to be from  $\sig^1_{c_i}$ to $\sig^1_{c_{i+1}}$ rather than the other way round. Now let $h':= \bigcup h_i$.

Let $R:=\{2^{-i} \mid \iin\}$.
We claim that 
\labtequ{uco}{$h' \restr (R \times \unin)$ is uniformly continuous.}
For suppose not. Then, there is some $\eps\in R_+$ and an infinite sequence of pairs $P_i=\{p_i, q_i\}$ of points in $R \times \unin$ such that $d( h'(p_i) , h'(q_i) ) > \eps$ \fe\ $i$ and the distance between $p_i$ and $q_i$ converges to 0. Note that \fe\ $ s \in R$ the subspace $\{ s \} \times \unin$ is compact, thus the function $h' \restr (\{ s \}) \times \unin$ is uniformly continuous by \Lr{HeCa}. This means that $\{ s \} \times \unin$ cannot contain an infinite subsequence of \seq{P}\ for any $ s \in R$. Even more, $\{ s \} \times \unin$ cannot meet an infinite subsequence of \seq{P}, because the distance between $p_i$ and $q_i$ converges to 0. It follows that $\{0\} \times \unin$ contains an accumulation point $(0,x)$ of \seq{P}, i.e.\ a point $(0,x)$ every neighbourhood of which contains infinitely many pairs $P_i$. 

Now let $\del$ be some (small) positive real number. Pick an $x' \in (\Q \cap \unin)$ \st\ $|x'-x|<\del/2$, and consider the open ball $O:=B_{\del}((0,x'))$ in $\unin \times \unin$. 

Let $R_O:= O \cap (R \times \{x'\})$. Choosing $\del$ small enough we can make sure that 
\labtequ{reps}{\fe\ $ s \in R_O$ there holds $\ell( \rho_s )<  r + \eps$,}
where $\rho_s: \unin \to X$ is defined by $x \mapsto h'( s ,x)$; indeed, $\rho_s$ coincides by definition with some $\sig^1_i$, and $\lim_i \ell(\sig^1_i) = r$.

As $O \ni x$, there is an infinite subsequence of \seq{P}\ contained in $O$.
Moreover, by \eqref{pwise} $h'(R_O)$ has a unique accumulation point in $X$. Thus we can find a pair $P_j=\{p_j,q_j\}$  in $O$ \st\ if $s$ (respectively, $s'$) is the element of $R$ for which $p_j \in \{ s \} \times \unin$ (resp., $q_j \in \{ s'\} \times \unin$) holds, then $d( h'( s ,x'), h'(s',x') ) < \eps/2$. 

Since $\rho_s$ coincides with some $\sig^1_i$, it has \cosp. As $||p_j , ( s ,x') || < 2\del$, this together with~\eqref{reps} implies $d( h'(( s ,x')), h'(p_j) ) < 2\del(r+\eps)$; similarly, we also have $d( h'(( s' ,x')), h'(q_j) ) < 2\del(r+\eps)$. Thus, by the triangle inequality applied to the four points $h'(p_j) ,  h'(( s ,x')), h'(( s' ,x'))$ and $h'(q_j)$ we obtain 
$$d( h'(p_j) , h'(q_j) ) \leq 
2\del(r+\eps) + \eps/2 + 2\del(r+\eps).$$
Since $\eps$ and $r$ are fixed and we can choose $\del$ freely, we can force this distance to be smaller than $\eps$ contradicting the choice of the $P_i$. This proves~\eqref{uco}.\comment{... and implies that the $\sig_{c_i}$ converge pointwise.} 

\medskip

The completion of $R \times \unin$ is $(R\cup \{0\}) \times \unin$; thus, by \eqref{univ} and \eqref{uco}, $h' \restr( R \times \unin )$ can be extended into a uniformly continuous function $h'': (R\cup \{0\}) \times \unin \to X$. Next, we prove that 
\labtequ{hcont}{$h:= h'\cup h''$ is continuous.}
Clearly, $h$ is continuous at any point in $(0,1] \times \unin$. So pick $x \in \{0\} \times \unin$ and $\eps\in \R_+$. By the continuity of $h''$, there is a basic open neighbourhood $O_\eps$ of $x$ in $R \times \unin$ that is mapped by $h''$ within the ball $B_{\eps/2}(h(x))$. Let $m_\eps\in\N$ be large enough that $h_i$ has width less that $\eps/2$ for every $i\geq m_\eps$; such an $m_\eps$ exists by the second sentence of \Lr{Uc} and the choice of the $h_i$. Assume \obda\ that $O_\eps$ does not meet ${2^{-i}} \times \unin$ for $i<m_\eps$. Extend $O_\eps$ into a set $O'\subseteq \unin \times \unin$ as follows. For every $i\geq m_\eps$ and every point $p=({2^{-i}},y)\in O$, put into $O'$ the line segment $L_p$ connecting $p$ to the point $ (2^{-(i+1)},y)$. Note that for every point $y\in L_p$ we have $d( h'(y), h'(p) ) \leq \eps/2$ since $h'$ coincides with $h_i$ on $L_p$ by the definition of $h'$ and $h_i$ has width less that $\eps/2$. As $O \cap (R \times \unin)$ is mapped by $h''$ within the ball $B_{\eps/2}(h(x))$, this implies that $h(O')\subseteq B_{\eps}(h(x))$.

But $O'$ contains by construction an open subset of $\unin \times \unin$ containing $x$. This proves \eqref{hcont}, which means that $h$ is a homotopy in $X'$ between the closed 1-simplex $h(0,x)$ and $\sig^1_{m_0}= h(1,x)$. We now define $z(x):= h(0,x)$, which is going to be the simplex we are looking for.

\medskip

Note that \fe\ $j$ the restriction $h \restr ([0,2^{-j}] \times \unin)$ is a homotopy between $z$ and  $\sig^1_{m_j}$ in $X'$, but this homotopy does not use the \arex s $X'_1,\ldots, X'_{j-1}$. Thus, as the \arex\ $X'_i$ has by construction \evol\ $2^{-i}$ \fe\ $i$, we obtain $d_1( \ec{\sig^1_{m_j}} , \ec{z} )\leq 2^{-(j-1)}$ \fe\ $j$ by the definition of $d_1$, since $\sig_{m_j}$ and $z$ are homotopic in the \arex\ $\bigcup_{i\geq j} X'_i$ of $X$. This proves that
\labtequ{Z}{$(\ec{\sig^1_{m_i}})_{i\in\N}$ is a \Cs\ with limit $Z:=\ec{z}$.}





\subsection*{Step III}

Our next aim is to prove that 
\labtequ{zr}{$\ell(z) \leq r$.}
Recall that $r$ was defined in Step I. 
Suppose, to the contrary, there is a finite sequence $S=s_1 < s_2 < \ldots < s_k$ of points in $\unin$ with $\sum_{1\leq i< k} d(z(s_i),z(s_{i+1})) =: r' >r$. Clearly, we may assume that $s_j\in \Q$ \fe\ $j$. Let $\eps:= \frac{r'-r}{2k}$. By \eqref{pwise} and the construction of $h$ we obtain that $\lim_i \sig^1_{\beta_i}(s_j) = h(0, s_j) = z(s_j)$ \fe\ $j$. Thus, choosing $i_0\in\N$ large enough we can make sure that $d( \sig^1_{\beta_i} (s_j) , z(s_j) ) < \eps$ \fe\ $j$ and every $i>i_0$. But then, the sequence $S$ witnesses the fact that $\ell(\sig^1_{\beta_i}) \geq r'$ \fe\ $i>i_0$, which contradicts the choice of $\seq{\sig^1}$ and proves~\eqref{zr}.

From \eqref{zr} we will now easily yield
\labtequ{Zr}{$\ell(Z) = r$.}
Firstly, note that by \eqref{Z} and the definition of $\ell(Z)$ we have $\ell(Z) \leq r$ by \eqref{zr}. Suppose that $\ell(Z)=r' < r$, and let $(\ec{\sig'_i})_\iin$ be a \Cs\ in $Z$ with $\lim \ell(\sig'_i)=r'$. 
Replacing $\sig^1_{c_i}$ in $\chi_{c_i}$ \fe\ $i$ by $\sig'_i$ we obtain a new sequence $\seq{\chi'}$ from $\seq{\chi}$, and it follows easily from \eqref{Z} that $(\ec{\chi'_i})_\iin \in D$ since $(\ec{\chi_i})_\iin \in D$. But $\lim_i \ell(\chi'_{c_i}) = \lim_i\ell(\chi_{c_i}) - r + r' < \lim_i\ell(\chi_i)$, which contradicts the choice of \seq{\chi}. Thus $\ell(Z) = r$ as claimed.

\medskip

Similarly to the proof of \eqref{zr} one can also easily prove that
\labtequ{zwp}{$z$ has \cosp.}

We now claim that $Z\hleq D$. Indeed, we have $\ell(D - Z) \geq \ell(D) - \ell(Z)$ by \Cr{ltriang}. Moreover, by \Lr{sum} we have $D - Z = \lim (\ec{\chi_{a_i} - \sig^1_{a_i}})$, and thus  $$\ell(D - Z) \leq \lim \ell(\chi_{a_i} - \sig^1_{a_i}) = \lim \ell(\chi_{a_i}) - \lim \ell(\sig^1_{a_i}) = \ell(D) - r = \ell(D) - \ell(Z),$$
where we used \eqref{Zr}. Thus $Z\hleq D$ as claimed, and as $D$ is primitive we obtain $Z=D$.

\medskip

Finally, we claim that $z$ is a \clex. Easily, the simplex $z$ is closed since all the $\sig^1_i$ are. Suppose the image of $z$ is not a circle. Then, there must be points $x\neq y \in [0,1)$ \st\ $z(x)=z(y)$. Now consider the two simplices $z_1$ and $z_2$ obtained by subdividing $z$ at these two points $x,y$,  and define $Z_1:= \ec{z_1}$ and $Z_2:= \ec{z_2}$. Easily, $\ell(z)= \ell(z_1) + \ell(z_2)$. We will show that $Z_1 \hleq Z$. For this, note that $Z - Z_1 = \ec{z - z_1} $ by \Lr{sum}, and so $Z-Z_1= \ec{z_2} = Z_2$. Thus 
$$\ell(Z - Z_1) = \ell(Z_2) \leq  \ell(z_2) = \ell(z) - \ell(z_1) \leq \ell(z) - \ell(Z_1) = \ell(Z) - \ell(Z_1),$$
and with \Cr{ltriang} we obtain $\ell(Z - Z_1) =\ell(Z) - \ell(Z_1)$, i.e.\ $Z_1\hleq Z$ as claimed. But as we have already shown that $Z=D$ and $D$ was assumed to be primitive, we obtain $Z = Z_1$, and thus $\ell(z_1) = \ell(z)$ since $\ell(z)=\ell(Z)$. This means that  $\ell(z_2) = 0$, which cannot be the case by~\eqref{zwp}. This contradiction proves that $z$ is a \clex.


\comment{
	\begin{lemma}\label{triang}
	Let $\seq{\sig}$ be an infinite sequence of \adim\ simplices in a compact space $X$ \st\ $0< \lim \ell(\sig_i)< \infty$. Then, \fe\ $\eps>0$ there is a $k\in\N$, a subsequence \susq{\sig}{a}\ of $\seq{\sig}$, and a triangulation of each $\sig_{a_i}$ into $k$ simplices $\sig_i^1, \ldots \sig_i^k$ with the following properties: 
\begin{itemize}
\item \fe\ $j\in \{1,\ldots,k\}$, if $s_i^j,t_i^j$ are the vertices of $\sig_i^j$ 
then the sequences \susq{s^j}{a}\ and \susq{t^j}{a}\ converge.
\item \fe\ $j\in \{1,\ldots,k\}$ and every $\iin$ there holds $\ell(\sig_{a_i}^j)<\eps$.
\item \fe\ $j\in \{1,\ldots,k\}$ there holds $0< \lim_i \ell(\sig_{a_i}^j) < \infty$.
\end{itemize}
	\end{lemma} 
}

\end{proof}

Thus we have proved \Lr{lem}, which combined with \Cr{split} proves our main result \Tr{lCycD}:

\begin{proof}[Proof of \Tr{lCycD}.] 
Suppose first that $\ell(C)< \infty$. Then we can apply \Cr{split} to obtain $C= \blcup_{i\in I} D_i$ where the $D_i$ are \prim. Applying \Lr{lem} to each $D_i$ we obtain a \clex\ $z_i$ with $D_i = \ec{z_i}$ and $\ell(z_i)=\ell(D_i)$. Note that we have $\ell(C)= \sum \ell(D_i)$ by the definition of $\blcup$. Thus $\ell(C)=\sum \ell(z_i)$. It remains to check that \seq{z}\ is a \rep\ of $C$. Indeed, we have $C = \lim \sum_{j\leq i} D_i$ by \eqref{ctblsum}, and substituting $D_i$ by $\ec{z_i}$ we obtain
$C = \lim \sum_{j\leq i} \ec{z_i}$, which means that \seq{z}\ is indeed a \rep\ of $C$ by definition. This proves the assertion in this case.

The other case, when $\ell(C)= \infty$ is easier. All we need to show is the existence of a \rep\ of $C$. For this, let \seq{C}\ with $C_i\in \homp$ be a sequence in $C$, and \fe\ $C_i$ pick an 1-cycle $c_i$ \st\ $\ec{z_i} \in C_i$. Now putting $z_i:= c_i - \sum_{j<i} c_j$, we obtain a \rep\ \seq{z}\ of $C$.
\end{proof}

\section{Application to graphs} \label{Sappg}
\newcommand{\lER}{\ensuremath{\ell: E(G) \to \R_{>0}}}
\newcommand{\ltopx}[1]{\ensuremath{ |#1|_\ell }}
\newcommand{\ltopxl}[2]{\ensuremath{|#2|_{#1} }}
\newcommand{\ltp}{\ltopx{$G$}}
\newcommand{\ltpf}[1]{\ltopf{#1}($G$)}

In this section we show that the \tcs\ \ccg\ described in the Introduction can be obtained as a special case of \lhom, and that our main result implies, in fact strengthens, \Tr{Zer}. The reader of this somewhat technical section is expected to be familiar with \ccg\ and the terminology and ideas of \cite[Chapter~8.5]{DiestelBook05}.

Let us first prove 

\begin{theorem} \label{lhiscc}
 For every \lfg\ \g \ti\ a metric of \fcg\ \st\ $\lhom(\fcg)$ is canonically isomorphic to \ccg.
\end{theorem}

The metric $d_\ell$ we are going to use in \Tr{lhiscc} is induced by an assignment $\lER$ of lengths to the edges of \G. More precisely, any such assignment naturally induces a distance $d_\ell(x,y)$ between any two points $x,y$, and we let \defi{\ltp} denote the completion of the corresponding metric space. For more details see \cite{ltop}, where the space \ltp\ is extensively studied. It turns out that choosing an appropriate assignment $\ell$ one obtains a metric space homeomorphic to \fcg:

\begin{theorem}[Georgakopoulos \cite{ltop}] \label{finlf}
If $G$ is \lf\ and\\ $\sum_{e \in E(G)} \ell(e) < \infty$ then $\ltp \cong \fcg$.
\end{theorem}

\begin{proof}[Proof of \Tr{lhiscc} (sketch)]
Fix a normal spanning tree $T$ of \G. Choose \lER\ \st\ $\ltp \cong \fcg$ and moreover the sums of the squares of the lengths of the fundamental cycles \wrt\ $T$ is finite. For example, we could start with an assignment $\ell'$ with $\sum \ell'(e)<\infty$, which guarantees $\ltp \cong \fcg$  by \Tr{finlf}, and then let $\ell(e):= \ell'(e)/m(e)$ where $m(e)$ is the number of fundamental cycles containing $e$.

We now define a map $f: \ccg \to \lhom(\ltp)$ which will turn out to be a canonical isomorphism.
Given a $C\in \ccg$,  write $C$ as the sum of a family \cf\ of fundamental cycles \wrt\ $T$; this is possible by \cite[Theorem~8.5.8]{DiestelBook05}. We will now construct a loop \sig\ in \ltp\ whose class will become the image $f(C)$ of $C$. We begin with a loop $\tau$ in \ltp\ that traverses each edge of $T$ once in each direction and traverses no other edges of \G. To see that such a loop exists, replace each edge of $T$ by a pair of parallel edges to obtain the auxiliary multigraph $T'$, and apply \cite[Theorem~2.5]{RDsBanffSurvey} to obtain a topological Euler tour $\tau'$ of $T'$. Now $\tau'$ clearly `projects' to the desired loop $\tau$. We then modify $\tau$ into \sig\ by attaching to it the cycles in \cf. To achieve this, assume that $\tau$ maps a non-trivial interval $I_v$ to each vertex $v$ of \G. Now for every fundamental cycle $F\in \cf$, let $v_F w_F$ be the chord of $F$, and assume without loss of generality that $v_F$ is closer to the root of $T$ than $w_F$. Modify $\tau$ so as to use the interval $I_{v_F}$, previously mapped to $v_F$, in order to travel once around $F$, starting and ending at $v_F$. Doing so for every $F\in \cf$ we obtain the loop \sig\ from $\tau$. One still has to check that \sig\ is indeed continuous, but this is not hard.
We let $f(C):= \ec{1\sig} \in \lhom(\fcg)$. 

The map $f$ is well-defined since $T$ and $\tau$ are fixed, and every $C\in \ccg$ has a unique representation as a sum of fundamental cycles \wrt\ $T$. 

To see that $f$ is injective, let $C \neq D \in \ccg$. Then the representations of $C$ and $D$ as sums of fundamental cycles differ by at least one fundamental cycle, since there must be a chord $e$ of $T$ contained in one of $C,D$ but not in the other. Now following the lines of \Tr{TSi} one can prove that $f(C)\neq f(D)$; indeed, $d_1(f(C),f(D))$ is bounded from below by a function of the length of $e$.

It remains to show that $f$ is onto. Pick an element $B$ of $\lhom(\fcg)$ for which we would like to find a preimage. Let \seq{B}\ be a \Cs\ in $B$. \Fe\ $B_i$ choose an 1-cycle $\chi_i$ \st\ $\ec{\chi_i}= B_i$. Using the loop $\tau$ from our earlier construction, we can join all the simplices in $\chi_i$ into one loop $\rho_i$ which, as $\tau$ is null-homotopic, is homologous to $\chi_i$. Now let $C_i \in \ccg$ be the sum $\sum \{a_e F_e \mid e\in E(G)\sm E(T)\}$ of fundamental cycles whose chords are traversed by $\rho_i$ (here $F_e$ denotes the fundamental cycle containing the chord $e$ and $a_e$ is the multiplicity of traversals of $e$ by $\rho_i$).

We claim that $f(C_i)$ is the equivalence class of the constant sequence $(\ec{\rho_i})$. To begin with, recall that $f(C_i)$ is by definition the equivalence class of the constant sequence $(\ec{\sig_i})$ for some loop $\sig_i$ that traverses the same chords of $T$ as $\chi_i$ does. However, the two loops will in general not be homologous, since the order in which these chords are traversed may differ at infinitely many positions. But \lhom\ has the ability of `disentangling' infinite products of commutators, and indeed, we will show that $d_1( \ec{\rho_i}, \ec{\sig_i} ) = 0$. For this, recall that we chose the edge-lengths $\ell(e)$ so that the sum of the squares of the lengths of all the fundamental cycles is finite. Applying \Lr{U} to each fundamental cycle, we can construct an \arex\ of \ltp\ with finite \evol\ in which every fundamental cycle is null-homologous. This means that \fe\ $\eps>0$ \ti\ an \arex\ $X_\eps$ of \ltp\ of \evol\ at most \eps\ in which all but finitely many fundamental cycles are null-homologous. Note that in each such $X_\eps$ the loops $\rho_i$ and $ \sig_i$ are homologous, since they traverse the same chords, and all but finitely many of these chords do not matter in $X_\eps$; thus the order in which they traverse the chords does not matter (recall that $H_1$ is abelian). This means that $d_1( \ec{\rho_i}, \ec{\sig_i} ) = 0$ as claimed.

We have thus found a sequence $C_i \in \ccg$ \st\ $(f(C_i))$ converges to $B$, but we would like to have an element $C\in \ccg$ with $f(C)= B$. To achieve this, we choose a subsequence $(C_{a_i})$ of $(C_i)$ that converges, as a set, to an element $C$ of \ccg; such a subsequence exists by compactness. It is now straightforward to check that $f(C)= B$ as desired: we can bound $d_1( f(C), f(C_{a_i}) )$ from above by any \eps\ choosing $i$ large enough. Indeed, choose $i$ so that the sum of the squares of the lengths of the fundamental cycles \wrt\ chords in the symmetric difference $ C - C_{a_i}$ is small compared to \eps. Since the sequence $(f(C_i))$ converges to $B$ this immediately yields $f(C)=B$. This completes the proof that $f$ is onto, which makes it an isomorphism, and by construction a canonical one.
\end{proof}

Using this, one now easily obtains \Tr{Zer} as a corollary of our main result \Tr{lCycD}. Indeed, given $C\in \ccg$ we apply \Tr{lCycD} to $f(C)$, where $f$ is the canonical isomorphism of \Tr{lhiscc}, to obtain a \rep\ $(z_i)$ of $f(C)$ with every $z_i$ being a \clex. Now if two of these \clex es share an edge $e$, then we can remove $e$ from both and combine the remaining arcs into a new closed simplex, thus obtaining  a new \rep\ of smaller total length, contradicting \Tr{lCycD}. This proves that the $z_i$ are edge-disjoint, and since $f$ is canonical the $f^{-1}(z_i)$ correspond to the same circles of \fcg\ and sum up to $C$.

In fact, this way we get something slightly stronger than \Tr{Zer}: for a given $C\in \ccg$ there may be several ways to decompose it as a sum of edge-disjoint circles; see \cite[p. 6]{hotchpotch} for an interesting example. \Tr{Zer} cannot distinguish between any of those ways, but our \Tr{lCycD} can: it returns one of minimal length. As the total length  of such a decomposition does not only depend on the edge-set (see \cite[Example 4.5.]{ltop}), this fact can be used in order to control the decomposition we obtain by varying the edge-lengths.

Furthermore, with \Tr{mainI} we generalise, in a sense, \Tr{Zer} to \nlfg s. For such graphs there are many candidate topologies on which \ccg\ can be based, so there is no standard cycle space theory. \Tr{mainI} helps to overcome this difficulty by offering a general result that, for each choice of a topology, yields a corollary similar to \Tr{Zer}. This approach is explained in \cite[Section 5]{ltop}.

\section{Higher dimensions} \label{Shidi}


Our definition of \lhom\ can be easily adapted to yield higher dimensional homology groups \lhomn. One can then ask if an analogue of our main result \Tr{lCycD} still holds in higher dimensions, but one should first choose a notion of $n$-dimensional content $vol()$, since there are several ways to generalise `length' to higher dimensions. Having chosen such a notion, e.g.\ the $n$-dimensional Hausdorff measure, one could then try to prove the following. 

\begin{problem}\label{lCycDn}
For every compact metric space $X$ and $C \in \lhomn(X)$, there is a \rep\ $\seq{z}$ of $C$ with $\sum_i vol(z_i)=vol(C)$. 
\end{problem} 

Most parts of our proof \Tr{lCycD}, in particular \Tr{Tmezzo}, could still be used in an attempt to prove \Prb{lCycDn}. To begin with, one would need to generalise the results of \Sr{Siso} for the chosen notion of content, which does not seem to be hard. The biggest difficulty though seems to be a generalisation of \eqref{uco}. 

\comment{
	difficulties:

\begin{itemize}
 \item need to either assume or prove an isoperimetric inequality: a cycle of small area can be patched by smaller volume
\item need to assume, or prove using (i), that two cycles of fixed area that are arbitrarily close can be made equivalent by small volume
\item uniform continuity of $h'$ is not clear, even in its restriction like in \eqref{uco}. 
\item (minor) a minimal cycle need not be the boundary of a sphere
\end{itemize}

\begin{lemma}\label{}
Let $\seq{\sig}$ be an infinite sequence of \ndim\ simplices in a compact space $X$. Then, \fe\ $\eps>0$ there is a subsequence \sseq{\sig}{a}\ of $\seq{\sig}$ and a triangulation of each $\sig_i$ into $k$ simplices $\sig_i^1, \ldots \sig_i^k$ with the following properties: 
\begin{itemize}
\item \fe\ $m\in \{1,\ldots,k\}$, and every $i,j\in \N$, if $u,v$ are corresponding vertices of $\sig_i^m$ and $\sig_j^m$ then $d(u,v)<\eps/kf_{1,n}$, and
\item \fe\ $m\in \{1,\ldots,k\}$, every $\iin$ and every face $F$ of dimension $n'$ of $\sig_i^m$ there holds $vol_{n'}(F)<\eps$.
\end{itemize}
\end{lemma} 
\begin{proof}
By induction on $n$
	\end{proof}
} 

\note{\section{} \begin{problem}\label{primsplit}
If $D$ is primitive, $D \hleq C$ for some $C$ and $C = C_1 \lcup C_2$ then $D\hleq C_i$ for some $i\in \{1,2\}$.
\end{problem} 

\begin{problem}
Inverse limits
\end{problem} 
}

\acknowledgement{I am very grateful to Reinhard Diestel for fruitful discussions that led to the definition of \lhom. I would like to thank the anonymous referee for their helpful remarks.}

\bibliographystyle{plain}
\bibliography{collective}
\end{document}